\documentclass[a4]{article}

\usepackage{amssymb}
\usepackage[pdftex]{graphicx}
\usepackage{epstopdf}
\usepackage{amsmath}
\usepackage{amsfonts}
\usepackage[english]{babel}
\usepackage{amsthm}
\usepackage{extarrows}
\usepackage{amscd}
\usepackage{tikz}
\usepackage{mathdots}
\usetikzlibrary{arrows,shapes}

\def\pf{\par\noindent {\bf Proof}~\par\noindent}
\def\qed{~\hfill{$\square$}\pagebreak[1]\par\medskip\par}

\newcommand{\mR}{\mathbb{R}}
\newcommand{\mC}{\mathbb{C}}

\newcommand{\mS}{\mathbb{S}}
\newcommand{\mH}{\mathbb{H}}
\newcommand{\mI}{\mathbb{I}}
\newcommand{\mJ}{\mathbb{J}}
\newcommand{\mK}{\mathbb{K}}
\newcommand{\mQ}{\mathbb{Q}}

\newcommand{\bG}{\boldsymbol{G}}

\newcommand{\bO}{\boldsymbol{O}}

\newcommand{\mcD}{\mathcal{D}}

\newcommand{\mcQ}{\mathcal{Q}}

\newcommand{\mcH}{\mathcal{H}}

\newcommand{\mcE}{\mathcal{E}}

\newcommand{\gf}{\mathfrak{f}}
\newcommand{\gfd}{\mathfrak{f}^{\dagger}}
\newcommand{\gsl}{\mathfrak{sl}}

\newcommand{\gosp}{\mathfrak{osp}}

\newcommand{\uX}{\underline{X}}
\newcommand{\uY}{\underline{Y}}
\newcommand{\uv}{\underline{v}}
\newcommand{\uvd}{\underline{v}^\dagger}
\newcommand{\ux}{\underline{x}}
\newcommand{\uXi}{\underline{\Xi}}
\newcommand{\uxi}{\underline{\xi}}

\newcommand{\puZd}{\partial_{\uzd}}
\newcommand{\bmDZ}{\boldsymbol{\mathcal{D}}_{(\uz,\uzd)}}

\newcommand{\bd}{\boldsymbol{\delta}}

\newcommand{\bdS}{\boldsymbol{d\Sigma}}
\newcommand{\uz}{\underline{z}}
\newcommand{\uzJ}{\underline{z}^J}
\newcommand{\uvJ}{\underline{v}^J}
\newcommand{\uzJd}{\underline{z}^{\dagger J}}
\newcommand{\uvJd}{\underline{v}^{\dagger J}}
\newcommand{\uzd}{\underline{z}^{\dagger}}

\newcommand{\upz}{\partial_{\uz}}
\newcommand{\upzJ}{\partial_{\uz}^J}

\newcommand{\upzd}{\partial_{\uz}^{\dagger}}
\newcommand{\upzJd}{\partial_{\uz}^{\dagger J}}

\newcommand{\p}{\partial}
\newcommand{\dirac}{\p}

\newcommand{\eop}{\hfill$\square$}

\newcommand{\SE}{\mcE}

\newcommand{\onehalf}{\frac{1}{2}}

\newtheorem{theorem}{Theorem}
\newtheorem{lemma}{Lemma}
\newtheorem{proposition}{Proposition}
\newtheorem{definition}{Definition}
\newtheorem{remark}{Remark}
\newtheorem{corollary}{Corollary}

\newtheorem{example}{Example}

\date{}

\title{The Cauchy Integral Formula  in Hermitian, Quaternionic\\ and $\gosp(4|2)$ Clifford analysis}
\small{
\author
{F.\ Brackx$^\ast$, H.\ De Schepper$^\ast$, R.\ L\'{a}vi\v{c}ka$^\dagger$, V.\ Sou\v{c}ek$^\dagger$}
\vskip 1truecm
\date{\small  $^\ast$ Department of Electronics and Information Systems\\ Faculty of Engineering
and Architecture, Ghent University\\ Technologiepark-Zwijnaarde 126,
B--9052 Ghent, Belgium \\[2mm]
$^\dagger$ Mathematical Institute, Faculty of Mathematics and Physics, \\Charles University, Sokolovsk\'{a} 83, 186 75 Praha, Czech Republic\\[5mm]}

\begin{document}

\maketitle

\begin{abstract}
As is the case for the theory of holomorphic functions in the complex plane, the Cauchy Integral Formula has proven to be a corner stone of Clifford analysis, the monogenic function theory in higher dimensional euclidean space. 
In recent years, several new branches of Clifford analysis have emerged. Similarly as hermitian Clifford analysis in euclidean space $\mR^{2n}$ of even dimension emerged as a refinement of euclidean Clifford analysis by  introducing  a complex structure on $\mR^{2n}$, quaternionic Clifford analysis arose as a further refinement by introducing a so--called hypercomplex structure $\mQ$, i.e.\ three complex structures ($\mI$, $\mJ$, $\mK$) which submit to the quaternionic multiplication rules, on $\mR^{4p}$, the dimension now being a fourfold. Two, respectively four, differential operators lead to first order systems invariant under the action of the respective symmetry groups U$(n)$ and Sp$(p)$. Their simultaneous null solutions are  called hermitian monogenic and quaternionic monogenic functions respectively. In this contribution we further elaborate on the Caychy Integral Formula for hermitian and quaternionic monogenic functions. Moreover we establish Caychy integral formul{\ae} for $\gosp(4|2)$--monogenic functions, the newest branch of Clifford analysis refining quaternionic monogenicity by taking  the underlying symplectic symmetry fully into account. \\

\end{abstract}


\section{Introduction}
\label{introduction}


The concept of fundamental solution of a differential operator is crucial to the development of the function theory for the null--solutions of this operator. In particular the Cauchy kernel
$$
E(z) = \frac{1}{2\pi \, i} \frac{1}{z}
$$
which is the fundamental solution of the Cauchy--Riemann operator $\partial_{\overline{z}} = \onehalf\, (\p_x + i\, \p_y)$, i.e.
$$
\partial_{\overline{z}} E(z) = \delta(z)
$$
is at the heart of the theory of the holomorphic functions in the complex plane. A corner stone in this theory is the Cauchy Integral Formula, which reproduces a holomorphic function $f$ in the interior of a bounded domain $D$ from its values on the  (piecewise) smooth boundary $\partial D$:
\begin{equation}
\label{cauchycomplex}
f(z)  =  \frac{1}{2\pi i}\, \int_{\partial  D} \frac{f(\xi)}{\xi-z}\, d\xi, \quad z \in \overset{\circ}{D}
\end{equation}
The Cauchy kernel also is the key ingredient of the Cauchy transform which generates a holomorphic function $H$ in the interior and the exterior of $D$ from a given smooth function $h$ on the boundary $\partial D$, through the integral
$$
H(z)  =  \frac{1}{2\pi i}\, \int_{\partial  D} \frac{h(\xi)}{\xi-z}\, d\xi, \quad z \notin \p D
$$
\\

The Cauchy Integral Formula in the complex plane has been generalised to the case of several complex variables in two ways. On the one hand taking a holomorphic kernel and integrating over the distinguished boundary $\partial_0 D=\prod_{j=1}^n \partial D_j$ of a polydisk $D=\prod_{j=1}^n D_j$ in $\mC^n$ leads to the representation formula
$$
f(z_1,\ldots,z_n)  =  \frac{1}{(2\pi i)^n}\, \int_{\partial_0 D} \frac{f(\xi_1,\ldots,\xi_n)}{(\xi_1-z_1) \cdots (\xi_n-z_n)}\, d\xi_1 \wedge \cdots \wedge d\xi_n\ , \quad z_j \in \overset{\circ}{D}_j
$$
On the other hand  integrating over the (piecewise) smooth boundary $\partial D$ of a bounded domain $D$ in $\mC^n$ in combination with the Martinelli--Bochner kernel, see e.g.\ \cite{kyt,kra}, which is no longer holomorphic but still harmonic, results into the formula
\begin{eqnarray}
\label{CIF_II}
f(\vec{z}) & = & \int_{\partial D} f(\vec{\xi})\, U(\vec{\xi},\vec{z})\ , \quad \vec{z} = (z_1, z_2, \ldots, z_n) \in \overset{\circ}{D}
\end{eqnarray}
with
\begin{equation}
\label{mabokernel}
U(\vec{\xi},\vec{z}) = \frac{(n-1)!}{(2\pi i)^n}\, \sum_{j=1}^n\, (-1)^{j-1}\, \frac{\xi_j^c-z_j^c}{\left|\xi-z\right|^{2n}}\, \left [ d\xi_j \right ]
\end{equation}
where
$$
\left [ d\xi_j \right ] = d\xi_1^c \wedge \cdots \wedge d\xi_{j-1}^c \wedge d\xi_{j+1}^c \wedge \cdots \wedge d\xi_n^c \wedge d\xi_1 \wedge \cdots \wedge d\xi_n
$$
and $\cdot^c$ denotes complex conjugation. For some historical background on formula (\ref{CIF_II}), which was obtained independently by Martinelli and Bochner, we refer to \cite{kra}. When $n=1$ it reduces to the traditional Cauchy Integral Formula (\ref{cauchycomplex}); for $n>1$, it establishes a connection between harmonic and complex analysis.\\[-2mm]

An alternative for generalising the Cauchy Integral Formula to higher dimension is offered by Clifford analysis, which originally studied the so--called monogenic functions, i.e.\ continuously differentiable functions defined in an open region of euclidean space $\mR^m$, taking their values in the Clifford algebra $\mR_{0,m}$, or subspaces thereof, and vanishing under the action of the Dirac operator 
$$
\dirac = \sum_{\alpha=1}^m e_{\alpha}\, \p_{X_{\alpha}}
$$
which corresponds, under Fourier duality, to the Clifford vector variable
$$
\uX = \sum_{\alpha=1}^m e_{\alpha}\, X_{\alpha}
$$
$\left ( e_\alpha \right )_{\alpha=1}^m$ being an orthornomal basis of $\mR^m$ underlying the construction of the Clifford algebra $\mR_{0,m}$. Monogenic functions are the natural higher dimensional counterparts of holomorphic functions in the complex plane. The Dirac operator factorizes the Laplacian: $\Delta_m = - \dirac^2$, and is invariant under the action of the $\mathrm{Spin}(m)$--group which doubly covers the $\mathrm{SO}(m)$--group, whence this framework is usually referred to as euclidean (or orthogonal) Clifford analysis.  Standard references in this respect are \cite{rood, groen, guerl, gilm}. \\

\noindent
In euclidean Clifford analysis the Cauchy Integral Formula for a monogenic function $f$ in an open neighbourhood of the closure of a bounded domain $D$ in $\mR^m$ with smooth boundary $\partial D$ reads
\begin{equation}
f(\uX) = \int_{\partial D} E(\uY-\uX)\, d\sigma_{\uY}\, f(\uY)\ , \quad \uX \in \overset{\circ}{D}
\label{CCE}
\end{equation}
where now the Cauchy kernel $E(\uX)$ is the fundamental solution of the Dirac operator. It is given by
$$
E(\uX) = \frac{1}{a_{m}}\, \frac{\overline{\uX}}{\left| \uX \right|^{m}}
$$
$a_{m}$ being the area of the unit sphere $S^{m-1}$ in $\mR^{m}$, $\bar{\cdot}$ denoting Clifford algebra conjugation and $d\sigma_{\uX}$ being a Clifford algebra valued differential form of order $(m-1)$ given by
$$
d\sigma_{\uX} = \sum_{j=1}^m  e_j\, (-1)^{j-1}\, \widehat{dX_j} 
$$
where the notation $\widehat{dX_j} $ means that $dX_j$ is omitted in the outer product of the differentials, i.e.
$$
 \widehat{dX_j}  = dX_1 \wedge \ldots \wedge dX_{j-1} \wedge dX_{j+1} \wedge \ldots \wedge dX_n, \qquad j=1,\ldots,m
$$
For a thorough study of the concept of fundamental solution in Clifford analysis we refer to \cite{mexico}.\\

The Cauchy Integral Formula (\ref{CCE}), which reproduces a monogenic function in the interior of the domain $D$ from its values on $\p D$, has been  a corner stone in the function theoretic development of euclidean Clifford analysis. The related Cauchy transform acting on smooth functions $h$ on $\p D$, generates monogenic functions in the interior $D^+ = \overset{\circ}{D}$ and the exterior $D^- = \mR^m \backslash \overline{D}$ of $D$  through 
\begin{equation}
\label{cauchytfmonogenic}
g(\uX) = \int_{\partial D} E(\uY-\uX)\, d\sigma_{\uY}\, h(\uY)\ , \quad \uX \in D^+ \cup\,  D^-
\end{equation}
with 
$$
\lim_{\uX \rightarrow \infty}\, g(\uX) = 0
$$
the non--tangential boundary values on $\p D$ of which being given by Plemelj-Sokhotzki type formul\ae, see \cite{fbhds} or \cite[Section 3.5]{guerl}. 
\\

This paper is devoted to establishing a Cauchy Integral Formula for so--called $\gosp(4|2)$--monogenic functions, the newest branch in Clifford analysis, meanwhile giving an overview of the attempts to establish Cauchy integral formul{\ae}  in hermitian and quaternionic Clifford analysis. The ingredients in any of these settings should thus be: a differential operator $\mathcal{D}$ and its fundamental solution $\mathcal{K}$ serving as a kernel for an integral transform which will reproduce null solutions of $\mathcal{D}$ in the interior of a bounded domain $D$ out of  their values on the boundary $\p D$ of that domain, and also will generate null solutions of  $\mathcal{D}$ in the interior and the exterior of $D$ out of given function values on $\p D$.


\section{Hermitian monogenicity}
\label{hermitianmonogenic}


The first refinement of monogenicity is so--called hermitian monogenicity, for which the setting is fixed as follows: take the dimension to be even: $m=2n$, rename the variables as
$$
(X_1,\ldots,X_{2n}) = (x_1,y_1,x_2,y_2,\ldots,x_n,y_n)
$$ 
and consider the standard complex structure $\mI_{2n}$, i.e. the complex linear real $\mbox{SO}(2n)$--matrix 
$$
\mI_{2n} = \mathrm{diag} \begin{pmatrix} \phantom{-} 0 & 1 \\ -1 & 0 \end{pmatrix}
$$
for which $\mI_{2n}^2 = - E_{2n}$, $E_{2n}$ denoting the identity matrix. We then define the rotated vector variable $\uX_\mI$ and the corresponding rotated Dirac operator $\dirac_\mI$ by
\begin{eqnarray*}
\uX_\mI &=& \mI_{2n}[\uX] = \sum_{k=1}^n (-y_k e_{2k-1} + x_k e_{2k}) \\
\dirac_\mI &=& \mI_{2n} [ \dirac] = \sum_{k=1}^n ( - \p_{y_k} e_{2k-1} + \p_{x_k} e_{2k})
\end{eqnarray*}
A differentiable function $F$ taking values in the complex Clifford algebra $\mC_{2n}$ then is called {\em hermitian monogenic} in some open region $\Omega$ of $\mR^{2n}$, if and only if in that region $F$ is a solution of the system 
\begin{equation}
\label{systemhermonogenic}
\{ \dirac F = 0 ,  \dirac_\mI F = 0 \}
\end{equation}
However, one can also introduce hermitian monogenicity, involving a complexification, by means of the projection operators $\pi^\pm = \pm \frac{1}{2} ({\bf 1} \pm i \, \mI_{2n})$. They produce the Witt basis vectors 
\begin{eqnarray*}
\gf_k &=& \pi^- [ e_{2k-1} ] \ = \  - \frac{1}{2} ({\bf 1} - i \, \mI_{2n}) [e_{2k-1}], \qquad k=1,\ldots,n \\
\gfd_k &=& \pi^+ [ e_{2k-1} ] \ = \ \phantom{-}  \frac{1}{2} ({\bf 1} + i \, \mI_{2n}) [e_{2k-1}], \qquad k=1,\ldots,n
\end{eqnarray*}
submitting to the properties
$$
\gf_j \gf_k + \gf_k \gf_j = 0, \quad \gfd_j \gfd_k + \gfd_k \gfd_j = 0, \quad \gf_j \gfd_k + \gfd_k \gf_j = \delta{jk}, \qquad j,k=1,\ldots,n
$$
which imply their isotropy. By means of this Witt basis, we define the hermitian vector variables
\begin{eqnarray*}
\uz = \pi^-[\uX] = - \frac{1}{2} ({\bf 1} - i \, \mI_{2n} ) [\uX]  &=& \sum_{k=1}^n (x_k + i y_k)\, \gf_k \ = \ \sum_{k=1}^n z_k \gf_k, \\
\uzd = \pi^+[\uX] =\phantom{-} \frac{1}{2} ({\bf 1} + i \, \mI_{2n} ) [\uX] &=& \sum_{k=1}^n (x_k - i y_k)\, \gfd_k \ = \  \sum_{k=1}^n z_k^c \gfd_k
\end{eqnarray*}
having introduced complex variables $(z_k, z_k^c)$ in $n$ respective complex planes. Correspondingly, the hermitian Dirac operators arise:
\begin{eqnarray*}
\upzd& = & \frac{1}{2} \pi^-[\dirac] \  = \ - \frac{1}{4} ( {\bf 1} - i \, \mI_{2n})[\dirac] \ =\   \sum_{k=1}^n \p_{z_k^c} \gf_k, \\
\upz & =  & \frac{1}{2} \pi^+[\dirac] \ = \ \phantom{-} \frac{1}{4} ( {\bf 1} + i \, \mI_{2n})[\dirac] \ = \  \sum_{k=1}^n \p_{z_k} \gfd_k
\end{eqnarray*}
It follows that for a function $F$ on $\mR^{2n} \cong \mC^n$ the hermitian monogenic system (\ref{systemhermonogenic}) is equivalent to the system
$$
\{ \upz F = 0 , \upzd F = 0 \}
$$
which can be shown to be invariant under the action of the unitary group U$(n)$. The basics of hermitian monogenicity theory can be found in e.g.\ \cite{cs1,cs2,herm,PSS}; for group theoretical aspects of this function theory we refer to \cite{eel2,DES}. \\[-2mm]

In the real approach to hermitian monogenicity we have the fundamental solutions
$$
E(\uX) = \frac{1}{a_{2n}}\, \frac{\overline{\uX}}{|\uX|^{2n}}, \quad E_\mI(\uX) = \frac{1}{a_{2n}}\, \frac{\overline{\uX}_\mI}{|\uX_\mI|^{2n}}
$$
for the operators $\dirac$ and $\dirac_{\mI}$ respectively, where now $a_{2n}$ denotes the area of the unit sphere $S^{2n-1}$ in $\mR^{2n}$. By projection they give rise to their hermitian counterparts, explicitly given by:
\begin{eqnarray*}
E(\uz) &=& 2 \pi^-[E(\uX)] \ = \ -  E(\uX) + i\, E_\mI (\uX) \ = \  \frac{2}{a_{2n}}\, \frac{\uz}{\left| \uz \right|^{2n}}\\
E^\dagger(\uz) &=& 2 \pi^+[E(\uX)] \ = \  \phantom{- } E(\uX) + i\, E_\mI(\uX)  \ = \ \frac{2}{a_{2n}}\, \frac{\uz^\dagger}{\left| \uz \right|^{2n}}
\end{eqnarray*}
whence
$$
E(\uX) = \onehalf\, (E^\dagger(\uz) - E(\uz))
$$
However, the latter turn out not to be fundamental solutions for the hermitian Dirac operators. Indeed, it holds, in distributional sense, that
\begin{eqnarray}
\label{diraczEz}
\upz E(\uz) &=& \frac{1}{2p}\, \beta\,  \delta(\uz) + \frac{2}{a_{4p}}\,  \beta\, \mbox{Fp} \frac{1}{|\uz|^{4p}} - \frac{2}{a_{4p}}\,  (2p)\, \mbox{Fp} \frac{\uzd \uz}{|\uz|^{4p+2}}\\
\label{diraczdEdz}
\upzd E^\dagger(\uz) &=& \frac{1}{2p}\,  (2p-\beta)\, \delta(\uz) + \frac{2}{a_{4p}}\,  (2p-\beta)\, \mbox{Fp} \frac{1}{|\uz|^{4p}} - \frac{2}{a_{4p}}\,  (2p)\, \mbox{Fp} \frac{\uz \uzd}{|\uz|^{4p+2}}
\end{eqnarray}
where $\beta = \sum_{k=1}^n\, \gfd_k\, \gf_k$ is a Clifford constant and $\mbox{Fp}$ stands for the {\em finite parts} distribution.
But introducing the particular circulant $(2 \times 2)$ matrices
$$
\bmDZ \ = \ \left( \begin{array}{cc} \partial_{\uz}^{\phantom{\dagger}} & \puZd\\ \puZd & \partial_{\uz^{\phantom{\dagger}}} \end{array} \right), \quad
{\bf{E}}(\uz) \ = \  \left( \begin{array}{cc} E(\uz) & E^\dagger(\uz)\\ E^\dagger(\uz) & E(\uz) \end{array} \right), \quad
 {\bd}(\uz) \ = \  \left( \begin{array}{cc} \delta(\uz) & 0\\ 0 & \delta(\uz) \end{array} \right)
$$
it was obtained in \cite{knock} that 
$$
\bmDZ {\bf{E}}(\uz) = {\bd}(\uz)
$$ 
whence the concept of a fundamental solution has to be reinterpreted for a matrix Dirac operator. Also observe that the matrix Dirac operator $\bmDZ$ still factorizes the Laplacian, since
$$
4\, \bmDZ \bmDZ^\dagger =  {\boldsymbol \Delta_{2n}}
$$ 
where $\boldsymbol \Delta_{2n}$ denotes the diagonal matrix with the Laplace operator in dimension $2n$ as the diagonal element.\\

Consequently, also the concept of hermitian monogenicity has to be reinterpreted: we say that a circulant matrix function
$$
\bG_2^1 \ = \  \left( \begin{array}{cc} g_1 & g_2\\ g_2 & g_1 \end{array} \right) 
$$
with continuously differentiable entries $g_1$ and $g_2$ defined in an open region $\Omega \subset \mR^{2n}$ and taking values in $\mC_{2n}$, is hermitian monogenic  if and only if it satisfies in $\Omega$ the system 
$$
\bmDZ \bG_2^1 = \bO
$$ 
where $\bO$ denotes the zero matrix, or, explicitly, 
$$
\left\{
\begin{array}{ll}
\upz\, g_1 + \upzd\, g_2 = 0\\
\upzd\, g_1 + \upz\, g_2 = 0
\end{array}
\right.
$$
Clearly if the functions $g_1$ and $g_2$ both are hermitian monogenic then the circulant matrix function $\bG_2^1 $ is hermitian monogenic, but the converse does not hold in general. However for the special case of a diagonal or anti--diagonal matrix function
$$
\bG_0 =  \  \left( \begin{array}{cc} g & 0\\ 0 & g \end{array} \right) \quad {\rm or} \quad \bG_0 =  \  \left( \begin{array}{cc} 0 & g\\ g & 0 \end{array} \right)
$$
i.e.\ when $g_1=g$ and $g_2=0$ or vice versa, the hermitian monogenicity of $\bG_0$ coincides with the hermitian monogenicity of $g$. \\[-2mm]


\section{The Cauchy Integral Formula in the hermitian case}
\label{cauchyhermitian}


In the actual dimension the classical Cauchy Integral Formula (\ref{CCE}) still reads
$$
f(\uX) = \int_{\partial D} E(\uY-\uX)\, d\sigma_{\uY}\, f(\uY)\ , \quad \uX \in \overset{\circ}{D}
$$
but now with $E(\uX)$ as given in the previous section and the differential form $d\sigma_{\uX}$ of order $(2n-1)$ given by
$$
d\sigma_{\uX} = \sum_{j=1}^n \left (  e_{2j-1} \, \widehat{dx_j} - e_{2j} \widehat{dy_j} \right ) 
$$
according to the new notations.\\

\noindent
A formal Cauchy Integral Formula for hermitian monogenic circulant matrix functions was first obtained in \cite{knock}. We recall the consecutive steps needed to arrive at this result. Introducing the notations
\begin{eqnarray}
\widehat{dz_j} &=& dz_1 \wedge dz_1^c \wedge \ldots \wedge dz_{j-1} \wedge dz_{j-1}^c \wedge dz_j^c \wedge  dz_{j+1} \wedge dz_{j+1}^c \wedge \ldots \wedge dz_n \wedge dz_n^c \label{hatdz}\\
\widehat{dz_j^c} &=& dz_1 \wedge dz_1^c \wedge \ldots \wedge dz_{j-1} \wedge dz_{j-1}^c \wedge dz_j \wedge  dz_{j+1} \wedge dz_{j+1}^c \wedge \ldots \wedge dz_n \wedge dz_n^c \label{hatdzdagger}
\end{eqnarray}
it is easily obtained that 
\begin{eqnarray*}
\widehat{dz_j} &=& 2^{n-1} (-i)^n \left [ \widehat{dx_j} + i \widehat{dy_j} \right ] \\
\widehat{dz_j^c} &=& 2^{n-1} (-i)^n \left [ \widehat{dx_j} - i \widehat{dy_j} \right ] 
\end{eqnarray*}
leading to the hermitian differential forms defined to be
$$
d\sigma_{\uz^{\phantom{\dagger}}}  =  \phantom{-} \sum_{j=1}^n \gf_j^\dagger \, \widehat{dz_j}, \qquad
d\sigma_{\uz^{\dagger}}  =  - \sum_{j=1}^n \gf_j \, \widehat{dz_j^c}
$$
which also may be obtained by projection:
\begin{eqnarray*}
d\sigma_{\uz^{\phantom{\dagger}}}  & = & (-i)^n 2^{n-1} \pi^-[d\sigma_{\uX}] \ = \  - \frac{1}{2}  (-i)^n 2^{n-1}  \left ( d\sigma_{\uX} - i \, d\sigma_{\uX_{\mI}} \right ) \\
d\sigma_{\uz^{\dagger}}  & = & (-i)^n 2^{n-1} \pi^+[d\sigma_{\uX}] \ = \   \phantom{-} \frac{1}{2}  (-i)^n 2^{n-1}  \left ( d\sigma_{\uX} + i \, d\sigma_{\uX_{\mI}} \right ) 
\end{eqnarray*}
whence 
$$
d\sigma_{\uX} = \frac{1}{2^{n-1}(-i)^n}\, (d\sigma_{\uzd} - d\sigma_{\uz})
$$

\noindent
Then, for a bounded domain $D \subset \mR^{2n}$ with smooth boundary $\partial D$ and a hermitian monogenic full circulant matrix function $\bG_2^1$ in an open neighbourhood of $\overline{D}$, the Cauchy Integral Formula reads
\begin{equation}
\label{CIFhermitianmatrix}
 \bG_2^1(\uX) = \frac{1}{(-2i)^n} \int_{\partial D} \boldsymbol{E} (\uv - \uz)\, \bdS_{(\uv,\uvd)} \, \bG_2^1(\uY) \quad ,  \quad \uX \in \overset{\circ}{D} 
\end{equation}
where $\uv$ is the hermitian vector variable corresponding to $\uY = \uvd - \uv$ running over  $\partial D$, and $\uz$ is the hermitian vector variable corresponding to $\uX = \uzd - \uz$ situated in the interior of $D$. The matrix differential form $\bdS_{(\uv,\uvd)} $ is given by  
$$
\bdS_{(\uv,\uvd)} \ = \ \left( \begin{array}{cc}  d\sigma_{\uv^{\phantom{\dagger}}} &  d\sigma_{\uv^{\dagger}} \\ d\sigma_{\uv^{\dagger}}  & d\sigma_{\uv^{\phantom{\dagger}}} \end{array} \right)
$$
The multiplicative constant appearing at the right hand side of formula (\ref{CIFhermitianmatrix}) originates from the re-ordering of $2n$ real variables into $n$ complex planes.\\[-2mm]

In the special case where  $\bG_2^1$ is taken to be the diagonal matrix function $\bG_0$, the above formula reduces to a genuine Cauchy Integral Formula for the hermitian monogenic function $g$, which explicitly reads
\begin{equation}
g(\uX) = \frac{1}{(-2i)^n} \int_{\partial D}  \left [ E(\uv-\uz) d\sigma_{\uv} + E^\dagger(\uv-\uz) d\sigma_{\uv^{\dagger}} \right ] \, g(\uY)  \quad , \quad \uX \in \overset{\circ}{D} 
\label{CCH}
\end{equation}
together with the additional integral identity
\begin{equation}
\int_{\partial D}  \left [ E(\uv-\uz) d\sigma_{\uv^\dagger} + E^\dagger(\uv-\uz) d\sigma_{\uv} \right ] \, g(\uY)  = 0 \quad , \quad \uX \in \overset{\circ}{D} 
\label{CCHnv}
\end{equation}
which thus should be fulfilled by every  function $g$ which is hermitian monogenic in an open neighbourhood of $\overline{D}$.\\[-2mm]

\begin{remark}
{\rm
In formul{\ae} like (\ref{CIFhermitianmatrix}), (\ref{CCH}) and (\ref{CCHnv}) we have used, next to each other, the cartesian variables $\uX$ and $\uY$ on the one hand, and the hermitian variables $\uz$ and $\uv$ on the other. They are linked by the transition  formulae $\uX = \uzd - \uz$ and $\uY = \uvd - \uv$ respectively. Because the variables $\uX$ and $\uY$ represent points located in $\overset{\circ}{D} $ and on $\p D$ respectively, we have kept them as function arguments emphasizing their geometric location, but it is, quite naturally, well understood that the functions involved should be expressed in terms of the vector variables $(\uz, \uzd)$ and $(\uv, \uvd)$ respectively.
In the forthcoming integral formul{\ae} next in this paper the same notation convention will be used.
}
\end{remark}\hfill\\[-4mm]

Now we will follow another approach and show how the formul{\ae} (\ref{CCH}) and (\ref{CCHnv}) may be directly derived from the euclidean Cauchy Integral Formula (\ref{CCE}). To that end we consider functions taking values in complex spinor space
$$
\mS = \mC_{2n} \, I \cong \mC_n \, I
$$
which is realized here by means of the primitive idempotent $I = I_1 \ldots I_n$, with 
$$
I_j = \gf_j \gf_j^{\dagger}, \qquad j=1,\ldots,n
$$
In \cite{cs2} it was shown that complex spinor space $\mS$, considered as a U$(n)$--module, decomposes as
\begin{equation}
\mS = \bigoplus_{r=0}^n \mS^r = \bigoplus_{r=0}^n (\mC\Lambda_n^\dagger)^{(r)} I
\label{homparts}
\end{equation}
into the U$(n)$--invariant and irreducible subspaces
$$
\mS^r = (\mC\Lambda_n^\dagger)^{(r)} I, \qquad r=0,\ldots,n
$$
each of them consisting of $r$-vectors from $\mC\Lambda_n^\dagger$ multiplied by the idempotent $I$, where $\mC\Lambda_n^\dagger$ is the Grassmann algebra generated by the Witt basis elements $\{ \gf_1^\dagger,\ldots,\gf_n^\dagger\}$. The spaces $\mS^r$ are also called the {\em homogeneous parts} of spinor space. Consequently, any spinor valued function $g$ decomposes as
$$
g = \sum_{r=0}^n g^{r}, \quad g^{r} : \mC^n \longrightarrow \mS^r, \qquad r = 0,\ldots,n
$$
in its so-called homogeneous components. It is worth observing that the action of the hermitian Dirac operators on a function $F^r$ taking values in a fixed homogeneous part $\mS^r$, will have the following effect:
\begin{eqnarray*}
\upz F^r &:& \mC^n \longrightarrow \mS^{r+1}\\[2mm]
\upzd F^r &:& \mC^n \longrightarrow \mS^{r-1}
\end{eqnarray*}
whence for such a function, the notions of monogenicity and hermitian monogenicity are equivalent. Indeed, seen the fact that 
$$
\dirac = 2 ( \upz - \upzd )
$$
hermitian monogenicity clearly implies monogenicity for any differentiable function. Moreover for each homogeneous component $g^r$ taking values in the homogeneous part $\mS^r$, we have seen above that $\upz g^r$ will be $\mS^{r+1}$ valued, while $\upzd g^r$ will be $\mS^{r-1}$ valued, whence $\dirac g^r = 0$ will force both terms to be zero separately.\\

\noindent
 A similar decomposition, followed by an analysis of the values, may now be applied to the Cauchy Integral Formula (\ref{CCE}). Indeed, as all building blocks of the hermitian framework are obtained,  up to constants, by projection, we may, conversely, decompose
$$
E(\uX) = \frac{1}{2} \left ( E^\dagger(\uz) - E(\uz) \right ) \qquad \mbox{and} \qquad d\sigma_{\uX} = \frac{i^n}{2^{n-1}} \left ( d\sigma_{\uz^\dagger} - d\sigma_{\uz} \right )
$$
Substituting these into (\ref{CCE}) yields, for each $r = 0, 1, 2,\ldots, n$,
$$
g^r(\uX) = \frac{1}{(-2i)^n} \int_{\partial D} \left (E^\dagger(\uv - \uz) - E(\uv - \uz) \right ) \left ( d\sigma_{\uv^\dagger} - d\sigma_{\uv} \right ) g^r(\uY)
$$
or still
$$
g^r(\uX) = \frac{1}{(-2i)^n} \left[  \int_{\partial D}  \left( E^\dagger(\uv - \uz) d\sigma_{\uv^\dagger} + E(\uv - \uz)  d\sigma_{\uv}  \right) g^r(\uY) \right.
$$ 
$$
\hspace*{55mm} \left. -  \int_{\partial D}    \left( E^\dagger(\uv - \uz) d\sigma_{\uv} + E(\uv - \uz)  d\sigma_{\uv^\dagger}  \right) g^r(\uY) \right] 
$$
Seen the definitions of $E(\uz)$, $E^\dagger(\uz)$, $d\sigma_{\uz}$ and $d\sigma_{\uz^\dagger}$, we will have
$$
\left ( E^\dagger(\uv - \uz)\, d\sigma_{\uv^\dagger} + E(\uv - \uz)\,  d\sigma_{\uv}  \right)\, g^r(\uY) \, : \, \mS^r \longrightarrow \mS^r
$$
while
$$
E^\dagger(\uv - \uz)\, d\sigma_{\uv}\,  g^r(\uY) \, : \, \mS^r \longrightarrow \mS^{r+2}  \qquad \mbox{and} \qquad E(\uv - \uz)\,  d\sigma_{\uv^\dagger}\,  g^r(\uY) \, : \, \mS^r \longrightarrow \mS^{r-2} 
$$
We thus directly obtain (\ref{CCH}) for each homogeneous component $g^r$:
$$
g^r(\uX) = \frac{1}{(-2i)^n}   \int_{\partial D}  \left( E^\dagger(\uv - \uz) d\sigma_{\uv^\dagger} + E(\uv - \uz)  d\sigma_{\uv}  \right) g^r(\uY) \quad , \quad \uX \in \overset{\circ}{D} 
$$
 while (\ref{CCHnv}) can be replaced by the even stronger result
$$
\int_{\partial D}  E(\uv-\uz) d\sigma_{\uv^\dagger}  \, g^r(\uY) \ \ = \ \  0 \ \  = \ \ \int_{\partial D}  E^\dagger(\uv-\uz) d\sigma_{\uv} \, g^r(\uY)  \quad , \quad \uX \in \overset{\circ}{D} 
$$
since both terms take values in different homogeneous parts. Note that the latter identities are not trivial. Indeed, as $\uX \in \overset{\circ}{D} $, the integral kernels are not differentiable in $\overset{\circ}{D} $ whence the Stokes Theorem may not be applied.\\[-2mm]

\noindent
This conclusion may be directly generalised to any spinor valued function $g$. It suffices to decompose such a function into its homogeneous parts and invoke the fact that $g$ is hermitian monogenic if and only if all its homogeneous parts $g^{r}$ are. We may thus write the above results separately for each component $g^{r}$ and by adding them up we obtain the following result.

\begin{proposition}
Let the spinor--valued function $g$ be hermitian monogenic in the open region $\Omega \subset \mR^{2n}$. Then for each bounded domain $D$ with smooth boundary $\p D$, such that $\overline{D} \subset \Omega$, it holds that
$$
g(\uX) = \frac{1}{(-2i)^n}   \int_{\partial D}  \left( E^\dagger(\uv - \uz) d\sigma_{\uv^\dagger} + E(\uv - \uz)  d\sigma_{\uv}  \right) g(\uY) \quad , \quad \uX \in \overset{\circ}{D} 
$$
together with 
the non--trivial integral identities
\begin{align}
\label{identity1} \int_{\partial D}  E(\uv-\uz) d\sigma_{\uv^\dagger}  \, g(\uY) &= 0 \quad , \quad \uX \in \overset{\circ}{D} \\ 
\label{identity2}\int_{\partial D}  E^\dagger(\uv-\uz) d\sigma_{\uv} \, g(\uY)  &=0 \quad , \quad \uX \in \overset{\circ}{D}
\end{align}
\end{proposition}

Additional identities are obtained through the action of the hermitian Dirac operators $\upz$ and $\upzd$ on formula (\ref{CCH}):
\begin{align*}
0 &=   \int_{\partial D}\,   \upz\, E(\uv-\uz)\,  d\sigma_{\uv} \, g(\uY) \quad , \quad \uX \in \overset{\circ}{D}   \\
0 &=   \int_{\partial D}\,   \upzd\, E^\dagger(\uv-\uz)\,  d\sigma_{\uvd} \, g(\uY) \quad , \quad \uX \in \overset{\circ}{D} 
\end{align*}
Putting, for $\uv - \uz  \neq 0$,
\begin{align*}
K(\uv - \uz) &=  -\, \upz\, E(\uv-\uz) = \upzd\, E^\dagger(\uv-\uz) \\
&= \frac{2}{a_{2n}}\, \frac{1}{|\uv-\uz|^{2n+2}}\, \left( \beta\, (\uv-\uz)(\uvd-\uzd) + (\beta-n)(\uvd-\uzd)(\uv-\uz) \right)
\end{align*}
we obtain the non--trivial integral identities
\begin{align*}
0 &=   \int_{\partial D}\,   K(\uv-\uz)\,  d\sigma_{\uv} \, g(\uY) \quad , \quad \uX \in \overset{\circ}{D}   \\
0 &=   \int_{\partial D}\,   K(\uv-\uz)\,  d\sigma_{\uvd} \, g(\uY) \quad , \quad \uX \in \overset{\circ}{D} 
\end{align*}
involving a hermitian monogenic integral kernel with a pointwise singularity in $\overset{\circ}{D}$.

\begin{remark}
\label{remarkholo}
{\rm
As mentioned in the introduction, in complex analysis  an alternative way of generalising the Cauchy Integral Formula to higher dimension is by means of the Martinelli--Bochner kernel, see e.g.\ \cite{kyt,kra}. One of the remarkable results of hermitian monogenic function theory is that, when the considered functions take their values in the $n$th homogeneous part $\mS^n$ of complex spinor space, hermitian monogenicity coincides with holomorphy in the complex variables ($z_1,\ldots,z_n$) and the above hermitian Cauchy Integral Formula reduces to the Martinelli-Bochner formula, in this way establishing a nice and interesting  connection between hermitian Clifford analysis and complex analysis in several variables, see also Example 1 in the next section.
}
\end{remark}


\section{The Cauchy transform in the hermitian case}
\label{cauchytfhermitian}


Given a smooth function $h$ on the smooth boundary $\p D$ of the bounded domain $D$ in $\mR^{2n}$, our aim is to generate, through the Cauchy transform, a hermitian monogenic function in its interior $D^+$ and its exterior $D^-$.
To that end we consider the integral
$$
\int_{\partial D} \boldsymbol{E} (\uv - \uz)\, \bdS_{(\uv,\uvd)} \, \boldsymbol{H}_0(\uY), \quad  \uX \in D^+ \cup D^-
$$
for a diagonal matrix function $\boldsymbol{H}_0$ with $h$ as diagonal entry.  This integral results into a circulant matrix function $\boldsymbol{G}^1_2$ in $D^+ \cup\, D^-$, with entries $g_1$ and $g_2$, given by
\begin{align*}
g_1(\uX) &=    \int_{\p D}\, E(\uv-\uz)\, d\sigma_{\uv}\, h(\uY) + E^{\dagger}(\uv-\uz)\, d\sigma_{\uvd}\, h(\uY) \\
g_2(\uX) &=    \int_{\p D}\, E(\uv-\uz)\, d\sigma_{\uvd}\, h(\uY) + E^{\dagger}(\uv-\uz)\, d\sigma_{\uv}\, h(\uY)
\end{align*}
Action by the matrix operator $\bmDZ$ learns that the functions $g_1$ and $g_2$ satisfy in $D^+ \cup\, D^-$ the system
$$
\left\{
\begin{array}{ll}
\upz\, g_1 + \upzd\, g_2 = 0\\
\upzd\, g_1 + \upz\, g_2 = 0
\end{array}
\right.
$$
Now we restrict ourselves to considerations about the interior $D^+$ of the bounded domain $D$, the results for the exterior $D^-$ being completely similar. 
Apparently there are two possibilities for generating a hermitian monogenic function in $D^+$. The first possibility consists in assuming that the boundary function $h$ satisfies the condition $g_2(\uX) = 0, \forall \uX \in D^+$, whence the function $g_1(\uX)$ will be hermitian monogenic in $D^+$. The second possibility consists in assuming that $h$ satisfies the condition $g_1(\uX) = 0, \forall \uX \in D^+$ which entails the hermitian monogenicity of $g_2(\uX)$. At first sight both possibilities are of equal value, but it will become apparent that the first option has to be preferred.\\

Let us consider the first possibility and assume that the boundary function $h$ satisfies the condition
$$
g_2(\ux) = \int_{\p D}\, E(\uv-\uz)\, d\sigma_{\uvd}\, h(\uY) + E^{\dagger}(\uv-\uz)\, d\sigma_{\uv}\, h(\uY) = 0 \quad , \quad \forall \  \uX \in D^+
$$
turning
$$
g_1(\uX) =    \int_{\p D}\, E(\uv-\uz)\, d\sigma_{\uv}\, h(\uY) + E^{\dagger}(\uv-\uz)\, d\sigma_{\uvd}\, h(\uY) 
$$
into a hermitian monogenic function in $D^+$.
Assuming now that the boundary function $h$ is spinor valued, we can decompose it into its homogeneous parts:
$$
h = \sum_{r=0}^n\, h^r \quad , \quad h^r : \p D \rightarrow \mS^r
$$
whence the function $g_1$ may be rewritten as
\begin{align*}
g_1(\uX) &=   \sum_{r=0}^n\, \int_{\p D}\, E(\uv-\uz)\, d\sigma_{\uv}\, h^r(\uY) + E^{\dagger}(\uv-\uz)\, d\sigma_{\uvd}\, h^r(\uY) \\
&=  \sum_{r=0}^n\, g_1^r(\uX)
\end{align*}
where,
for each $r=0, \ldots, n$, the function $g_1^r$ 
takes its values in $\mS^r$and inherits its hermitian monogenicity in $D^+$ from $g_1$. In this way the Cauchy transform generates, out of the spinor--valued boundary function $h$, a hermitian monogenic function in  $D^+$, this construction being carried out {\em componentwise}.\\

Now consider the case where the spinor--valued boundary function $h$ satisfies the second condition
$$
\int_{\p D}\, E(\uv-\uz)\, d\sigma_{\uv}\, h(\uY) + E^{\dagger}(\uv-\uz)\, d\sigma_{\uvd}\, h(\uY) = 0
$$
Clearly this condition is fulfilled by the function $h$ if each of its homogeneous components $h^r$ does. But, for each $r=0,\ldots,n$, the hermitian monogenic function
$$
g_2^r (\uX) = \int_{\p D}\, E(\uv-\uz)\, d\sigma_{\uvd}\, h^r(\uY) + E^{\dagger}(\uv-\uz)\, d\sigma_{\uv}\, h^r(\uY)
$$
takes values in $\mS^{r-2}  \bigoplus \mS^{r+2}$, and the hermitian monogenic function $g_2 = \sum_{r=0}^n\, g_2^r$ in $\overset{\circ}{D}$ is not constructed componentwise. This explains our preference for the first approach.\\

Note that the same results may be obtained, starting from the euclidean Cauchy transform (\ref{cauchytfmonogenic}) for monogenic functions, by decomposing the Cauchy kernels and the differential forms into their hermitian counterparts. Indeed, let us consider a smooth function $h$ on $\partial D$ and let us assume from the beginning that it satisfies the condition
\begin{equation}
\label{firstcondition}
g_2(\uX) = \int_{\p D}\, E(\uv-\uz)\, d\sigma_{\uvd}\, h(\uY) + E^{\dagger}(\uv-\uz)\, d\sigma_{\uv}\, h(\uY) = 0 \quad , \quad \forall \  \uX \in D^+
\end{equation}
Its Cauchy transform, see (\ref{cauchytfmonogenic}), yields a monogenic function  in $D^+$:
\begin{equation}
\label{monogenicinterior}
g(\uX) = \int_{\p D}\, E(\uY-\uX)\, d\sigma_{\uY}\, h(\uY) \quad , \quad   \uX \in D^+
\end{equation}
Due to condition (\ref{firstcondition}) this function $g$ takes the form
$$
g(\uX) = \frac{1}{(-2i)^n}\,  \int_{\p D}\, E(\uv-\uz)\, d\sigma_{\uv}\, h(\uY) + E^{\dagger}(\uv-\uz)\, d\sigma_{\uvd}\, h(\uY) \quad , \quad   \uX \in D^+
$$
in which we recognize, up to a constant, the function $g_1$, in other words: the function $g$ is not only monogenic but also hermitian monogenic in $D^+$.\\[-2mm]

\noindent
It is interesting to investigate the boundary value of this hermitian monogenic function $g(\uX)$ for $\uX \in D^+$ tending non--tangentially to a certain point $\uXi$ on $\partial D$. Boundary values of monogenic functions have been studied in e.g.\ \cite{fbhds}; for a matricial treatment of the boundary behaviour of hermitian monogenic functions we refer to \cite{fbbdkhds}. Recall that the Cauchy transform of the smooth function $h$ on the boundary $\p D$ of the bounded domain $D$, given by
$$
g(\uX) = \int_{\p D}\, E(\uY - \uX)\, d \sigma_{\uY}\, h(\uY) \quad , \quad \uX \in D^+ \cup D^-
$$
belongs to the Hardy spaces $H_2(D^+)$ and $H_2(D^-)$ with non--tangential boundary values belonging to the Hardy spaces $H_2^+(\p D)$ and  $H_2^-(\p D)$ respectively, given by the Plemelj-Sokhotzki formul\ae
$$
\widetilde{g}^+(\uXi) = \lim_{D^+  \ni\, \uX \longrightarrow \Xi}\, g(\uX) = \phantom{-}\, \onehalf\, h(\uXi) + \onehalf\, \mH[h](\uXi) \quad , \quad \uXi \in \p D
$$
and 
$$
\widetilde{g}^-(\uXi) = \lim_{D^-  \ni\, \uX \longrightarrow \Xi}\, g(\uX) = -\, \onehalf\, h(\uXi) + \onehalf\, \mH[h](\uXi) \quad , \quad \uXi \in \p D
$$
where $\mH$ stands for the {\em Hilbert transform} given by
\begin{align}
 \mH[h](\uXi) &= 2\; {\rm Pv} \int_{\p D}\, E(\uY-\uXi)\, d\sigma_{\uY}\, h(\uY){\nonumber}\\
\label{hilbertdef} &= 2\, \lim_{\epsilon \rightarrow +0}\, \int_{\p D_\epsilon}\, E(\uY-\uXi)\, d\sigma_{\uY}\, h(\uY)
\end{align}
with $\p D_\epsilon = \{ \uY \in \p D : d(\uY, \uXi) > \epsilon  \}$.\\[1mm]
Now, as already mentioned above, for a boundary function $h$ satisfying the condition
\begin{equation*}
 \int_{\p D}\, E(\uv-\uz)\, d\sigma_{\uvd}\, h(\uY) + E^{\dagger}(\uv-\uz)\, d\sigma_{\uv}\, h(\uY) = 0 \quad , \quad \forall \  \uX \in D^+ \cup D^-
\end{equation*}
the Cauchy transform $g(\uX)$ becomes hermitian monogenic in $D^+ \cup D^-$ and so belongs to the Hardy spaces $\mcH_2(D^+)$ and $\mcH_2(D^-)$ of hermitian monogenic functions in $D^+$ and $D^-$ respectively, showing non--tangential boundary values in $L_2(\p D)$, see \cite{fbhds, fbbdkhds}. Introducing the Hardy space $\mcH_2(\p D)$ as the closure in $L_2(\p D)$ of the non--tangential boundary values of the functions in $\mcH_2(D^+)$, these considerations lead to the following result.

\begin{proposition}
\label{hardyspace}
If the function $h$ belongs to the Hardy space $H_2(\p D)$ and satisfies the integral condition
\begin{equation}
\label{hardycondition}
 \int_{\p D}\, E(\uv-\uz)\, d\sigma_{\uvd}\, h(\uY) + E^{\dagger}(\uv-\uz)\, d\sigma_{\uv}\, h(\uY) = 0 \quad , \quad \forall \  \uX \in D^+ \cup D^-
\end{equation}
then $h$ belongs to the Hardy space $\mcH_2(\p D)$
\end{proposition}

\noindent
As to the the converse of Proposition \ref{hardyspace} we are able to prove the following.

\begin{proposition}
\label{hardyspaceconverse}
If the spinor--valued function $h$ belongs to the Hardy space $\mcH_2(\p D)$, then $h$ trivially belongs to the Hardy space $H_2(\p D)$ and moreover satisfies the integral condition (\ref{hardycondition}).
\end{proposition}

\pf 
If a function $f$ belongs to the Hardy space $\mcH_2(\p D)$ then there exists a function $G \in \mcH_2(D^+) \subset H_2(D^+)$ such that 
$$
\lim_{\uX \rightarrow \uXi}\, G(\uX) = f(\uXi) \quad , \quad \uXi \in \p D
$$
On the other hand the Cauchy transform $F$ of $f$ belongs to $H_2(D^+)$ and shows the non--tangential boundary value $f$. This means that the function $F - G$ is monogenic in $D^+$ and shows the non--tangential boundary value $0$, whence $F = G$. In other words: the Cauchy transform of a function $f \in \mcH_2(\p D)$ is hermitian monogenic in $D^+$ and belongs to $\mcH_2(D^+)$ with non--tangential boundary value $f$. When applying this property to the $r$--homogeneous component $h^r$ of $h$, we find that the Cauchy transform of $h^r$, which a priori takes values in $\mS^r \oplus \mS^{r+2} \oplus \mS^{r-2}$, shows the non--tangential boundary value $h^r$ which takes values in only $\mS^r$. This implies that the Cauchy transform of $h^r$ also takes its values in only $\mS^r$, whence it is the $r$--homogeneous component of the Cauchy transform of $h$. It follows that each of the homogeneous components of the boundary function $h$, and hence also $h$ itself, satisfies the integral condition (\ref{hardycondition}) and even the stronger conditions (\ref{identity1}) and (\ref{identity2}).
\eop

\begin{remark}
{\rm
Proposition \ref{hardyspaceconverse} was formulated for spinor valued functions, however without loss of generality since  the complex Clifford algebra $\mC_{2n}$ decomposes into a direct sum of copies of spinor space $\mS$.
}
\end{remark}

In addition condition (\ref{hardycondition}) makes it possible to rewrite and simplify
expression (\ref{hilbertdef}) defining the Hilbert transform, in terms of the hermitian counterparts to the Cauchy kernel functions and the differential forms involved. To that end we need the results of the following lemmata.

\begin{lemma}
\label{specificintint}
Let $\uX$ be a point in the interior of a bounded domain $D$ with smooth boundary $\p D$. One has
\begin{itemize}
\item[(i)] $\frac{1}{(-2i)^n}\, \int_{\p D}\, E^\dagger(\uv-\uz)\, d\sigma_{\uvd}  + E(\uv-\uz)\, d\sigma_{\uv} = 1$
\item[(ii)] $\int_{\p D}\, E^\dagger(\uv-\uz)\, d\sigma_{\uv} = 0$
\item[(iii)] $\int_{\p D}\, E(\uv-\uz)\, d\sigma_{\uvd} = 0$
\end{itemize}
\end{lemma}
\pf
These formul{\ae} readily follow from the well--known result, see e.g. \cite{fbhds},
$$
\int_{\p D}\, E(\uY-\uX)\, d\sigma_{\uY} = 1 \quad , \quad   \uX \in \overset{\circ}{D}
$$\eop

\begin{lemma}
\label{specificintboundary}
Let $\uXi$ be a boundary point of a bounded domain $D$ with smooth boundary $\p D$. One has
\begin{itemize}
\item[(i)] $\frac{1}{(-2i)^n}\, {\rm Pv}\, \int_{\p D}\, E^\dagger(\uv-\uxi)\, d\sigma_{\uvd}  + E(\uv-\uxi)\, d\sigma_{\uv} = \onehalf$
\item[(ii)] $ {\rm Pv}\, \int_{\p D}\, E^\dagger(\uv-\uxi)\, d\sigma_{\uv} = 0$
\item[(iii)] $ {\rm Pv}\, \int_{\p D}\, E(\uv-\uxi)\, d\sigma_{\uvd} = 0$
\end{itemize}
\end{lemma}
\pf
These formul{\ae} readily follow from the well--known result, see e.g. \cite{fbhds},
$$
{\rm Pv}\, \int_{\p D}\, E(\uY-\uXi)\, d\sigma_{\uY} = \onehalf
$$\eop

\begin{lemma}
\label{pvzero}
If the smooth function $h$ defined on the boundary $\p D$ of a bounded domain $D$ satisfies condition (\ref{firstcondition}), viz.
$$
 \int_{\p D}\, E(\uv-\uz)\, d\sigma_{\uvd}\, h(\uY) + E^{\dagger}(\uv-\uz)\, d\sigma_{\uv}\, h(\uY) = 0 \quad , \quad \forall \  \uX \in D^+
$$
then, for each point $(\uxi, \uxi^\dagger) \equiv \uXi \in \p D$,  it holds that 
$$
{\rm Pv}  \int_{\p D}\, E(\uv-\uxi)\, d\sigma_{\uvd}\, h(\uY) + E^{\dagger}(\uv-\uxi)\, d\sigma_{\uv}\, h(\uY) = 0
$$
\end{lemma}
\pf
In view of Lemma \ref{specificintint} {\em (ii)} and {\em (iii)}, it follows from (\ref{firstcondition}) that
$$
\lim_{\uX \rightarrow \uXi}\, \int_{\p D}\, \left( E(\uv-\uz)\, d\sigma_{\uvd} + E^{\dagger}(\uv-\uz)\, d\sigma_{\uv}  \right)\, \left( h(\uY) - h(\uXi) \right) = 0 
$$
whence
$$
\lim_{\epsilon \rightarrow +0}\, \int_{\p D_{\epsilon}}\, \left( E(\uv-\uxi)\, d\sigma_{\uvd} + E^{\dagger}(\uv-\uxi)\, d\sigma_{\uv}  \right)\, \left( h(\uY) - h(\uXi) \right) = 0 
$$
Invoking Lemma \ref{specificintboundary} {\em (ii)} and {\em (iii)}, we thus obtain
$$
\lim_{\epsilon \rightarrow +0}\, \int_{\p D_{\epsilon}}\, \left( E(\uv-\uxi)\, d\sigma_{\uvd} + E^{\dagger}(\uv-\uxi)\, d\sigma_{\uv}  \right)\,  h(\uY) = 0 
$$\eop

\noindent
The results of the above lemmata now lead to the following.

\begin{proposition}
\label{hilberttransformcondensed}
For a smooth function $h$ defined on the smooth boundary $\p D$ of a bounded domain $D$, satisfying condition (\ref{hardycondition}), viz.
$$
 \int_{\p D}\, E(\uv-\uz)\, d\sigma_{\uvd}\, h(\uY) + E^{\dagger}(\uv-\uz)\, d\sigma_{\uv}\, h(\uY) = 0 \quad , \quad \forall \  \uX \in D^+ \cup D^-
$$
the expression for its Hilbert transform reduces to
$$
 \mH[h](\uXi) = \frac{2}{(-2i)^n}\, {\rm Pv}\, \int_{\p D}\, E(\uv-\uxi)\, d\sigma_{\uv}\, h(\uY) + E^{\dagger}(\uv-\uxi)\, d\sigma_{\uvd}\, h(\uY)
$$
\end{proposition}\vspace{3mm}

\noindent
Note that we may now check the well--known Hilbert transform $\mH[1] = 1$. The constant function $1$ indeed satifies condition (\ref{hardycondition}) due to Lemma \ref{specificintint} {\em (ii)} and {\em (iii)}. Next, invoking Lemma \ref{specificintboundary} {\em (i)}, we find
$$
 \mH[1](\uXi) = \frac{2}{(-2i)^n}\, {\rm Pv}\, \int_{\p D}\, E(\uv-\uxi)\, d\sigma_{\uv}\,  + E^{\dagger}(\uv-\uxi)\, d\sigma_{\uvd}\, = 2\, \onehalf = 1
$$\hfill\\

There is, however, still a third approach possible to establish a Cauchy transform in the hermitian monogenic setting. Indeed, if the smooth boundary function $h$ on $\p D$ is assumed to satisfy the condition
\begin{equation}
\label{hardyconditionbis} \int_{\p D}\, K(\uv-\uz)\, d \sigma_{\uY}\, h(\uY) = 0
\end{equation}
where, recall, the hermitian monogenic kernel $K$ is given by
$$
K(\uv-\uz) = -\, \p_{\uz}\, E(\uv-\uz) = \p_{\uzd}\, E^\dagger(\uv-\uz)
$$
then, for the monogenic Cauchy transform of $h$, viz.
\begin{align*}
g(\uX) &=  \int_{\p D}\, E(\uY-\uX)\, d \sigma_{\uY}\, h(\uY) \\
&= \onehalf\, \int_{\p D}\, \left( E^\dagger(\uv-\uz) - E(\uv-\uz) \right)\,  d \sigma_{\uY}\, h(\uY) \quad , \quad \uX \in D^+ \cup D^-
\end{align*}
it holds that
$$
\p_{\uzd}\, g(\uX) = \onehalf\, \int_{\p D}\, \p_{\uzd}\,  E^\dagger(\uv-\uz) \,  d \sigma_{\uY}\, h(\uY) = \onehalf\, \int_{\p D}\, K(\uv-\uz) \,  d \sigma_{\uY}\, h(\uY) = 0 \quad , \quad \uX \in D^+ \cup D^-
$$
and also $\p_{\uz}\, g(\uX) = 0$ in $D^+ \cup D^-$ since $\dirac\, g(\uX) = 2(\p_{\uz} - \p_{\uzd})\, g(\uX) = 0$, in other words: $g(\uX)$ is hermitian monogenic in $D^+ \cup D^-$.\\

\noindent
Condition (\ref{hardyconditionbis}) now is the key to the following alternative characterization of the Hardy space $\mcH_2(\p D)$ in terms of the Hardy space $H_2(\p D)$.

\begin{proposition}
A function $f$ belongs to the Hardy space  $\mcH_2(\p D)$ if and only if $f$ belongs to the Hardy space $H_2(\p D)$ and moreover satisfies the integral condition
\begin{equation}
\label{hardyconditionprop}
 \int_{\p D}\, K(\uv-\uz)\, d \sigma_{\uY}\, f(\uY) = 0 \quad , \quad \uX \in D^+
\end{equation}
\end{proposition}

\pf\\[-3mm]
(i) Suppose that $f \in H_2(\p D)$. Then its Cauchy transform $F$ belongs to $H_2(D^+)$ and shows the non--tangential boundary value $f$ on $\p D$. Now if $f$ satisfies condition (\ref{hardyconditionprop}) then, as was shown above, $F$ is hermitian monogenic in $D^+$, in other words: $F \in \mcH_2(D^+)$ and so $f \in \mcH_2(\p D)$.\\
(ii) Suppose that $f \in \mcH_2(\p D)$. The following reasoning is similar to the proof of Proposition \ref{hardyspaceconverse}. There exists a function $G  \in \mcH_2(D^+)$ the non--tangential boundary value on $\p D$ of which is precisely the function $f$. As $\mcH_2(\p D) \subset H_2(\p D)$ the Cauchy transform $F$ of $f$ belongs to $H_2(D^+)$ and shows the non--tangential boundary value $f$ on $\p D$. Consider the function $F - G$; it is monogenic in $D^+$ and shows the non--tagential boundary value $0$ on $\p D$, whence $F = G$ in $D^+$. So it holds  that $F$ is hermitian monogenic in $D^+$, whence for $\uX \in D^+$,
$$
\p_{\uzd}\, F(\uX) =  \onehalf\, \p_{\uzd}\, \int_{\p D}\, \left( E^\dagger(\uv-\uz) - E(\uv-\uz) \right)\,  d \sigma_{\uY}\, f(\uY) = \onehalf \int_{\p D}\, K(\uv-\uz)\, d \sigma_{\uY}\, f(\uY) = 0 
$$
and condition (\ref{hardyconditionprop}) follows.
\eop\\

We conclude this section by the following  nice example.

\begin{example}
{\rm
Consider the function 
$$
f = F\, \gfd_1 \gfd_2 \cdots \gfd_n\, I
$$
with F a scalar function defined in an open neighbourhood $\Omega$ of $\overline{D}$, $D$
being a bounded domain in $\mR^{2n}$ with a smooth boundary $\p D$. Clearly $f$ takes values in the $n$-th homogeneous part $\mS^n$ of spinor space.
As was already mentioned in Remark \ref{remarkholo}, it was proven in \cite{knock} that $f$ is hermitian monogenic in $\Omega$ if and only if $F$ is a holomorphic function in the complex variables $(z_1, \ldots,z_n)$ and does not depend on the complex conjugates $(z_1^c, \ldots, z_n^c)$. Assuming now that the function $F(z_1,\ldots,z_n)$ is indeed holomorphic in $\Omega$, or, equivalently, that the $\mS^n$-- valued function $f$ is hermitian monogenic in $\Omega$, the Cauchy Integral Formula (\ref{CCH}) and the additional integral identities (\ref{identity1}) and (\ref{identity2}) hold for this function $f$. Seen the algebraic structure of $f$, implying that $d \sigma_{\uv}\, f(\uY) = 0$, these formul{\ae} reduce to
\begin{equation}
\label{exampleCIF}
f(\uX) = \frac{1}{(-2i)^n}\, \int_{\p D}\, E^\dagger(\uv-\uz)\, d \sigma_{\uvd}\, f(\uY) \quad , \quad \uX \in D^+
\end{equation}
and 
\begin{equation}
\label{exampleidentity}
\int_{\p D}\, E(\uv-\uz)\, d \sigma_{\uvd}\, f(\uY) = 0 \quad , \quad \uX \in D^+
\end{equation}
which, as we know, lead to the Martinelli--Bochner representation formula (\ref{CIF_II}), viz.
$$
F(z_1,\ldots,z_n) = \int_{\p D}\, U(\vec{\xi}, \vec{z})\, F(\xi_1,\ldots,\xi_n
$$
and the additional identities
\begin{equation*}
\int_{\p D}\, F(z_1, \ldots,z_n)\, \frac{1}{\rho^{2n}}\, (v_j - z_j)\, \widehat{d v_k^c} = \int_{\p D}\, F(z_1, \ldots,z_n)\, \frac{1}{\rho^{2n}}\, (v_k - z_k)\, \widehat{d v_j^c}, \quad , \quad j \neq k
\end{equation*}

\noindent
Now consider the continuous boundary function $h^n = f|_{\p D}: \p D \rightarrow \mS^n$, for which it holds, invoking (\ref{exampleidentity}), that 
$$
\int_{\p D}\, E(\uv-\uz)\, d \sigma_{\uvd}\, h^n(\uY) = 0 \quad , \quad \uX \in D^+
$$
meaning that the boundary function $h^n$ satisfies the first assumption.
It follows that the function $g$ given by
\begin{equation}
\label{exampleg}
g(\uX) = \frac{1}{(-2i)^n}\, \int_{\p D}\, E^{\dagger}(\uv - \uz)\, d \sigma_{\uvd}\, h^n(\uY) 
\end{equation}
is hermitian monogenic in $D^+$. Putting
$
g(\uX) = G\, \gfd_1 \cdots \gfd_n\, I 
$, 
the function $G$ then is holomorphic in the same region. 
Moreover it is clear that, seen (\ref{exampleCIF}), the function $g$ coincides with the initial hermitian monogenic  function $f$ in $D^+$,  while $g = 0$ in $D^-$.
It follows that, for $\uX  \in D^+ \longrightarrow \uXi \in \p D$, 
$$
\lim_{\uX \rightarrow \uXi}\, g(\uX) =  \lim_{\uX \rightarrow \uXi}\, f(\uX) 
$$
or
$$
 \onehalf\,  h^n(\uXi) + \onehalf\,  \mH[h^n](\uXi) = h^n(\uXi)
$$
or still
$$
\mH[h^n](\uXi) = h^n(\uXi)
$$
confirming that the restriction $h^n$ to $\p D$ of the hermitian monogenic function $F$ in $\Omega \supset \overline{D}$ belongs to the Hardy space $H_2^+(\p D)$, see \cite{fbhds, fbbdkhds}. Note that for the non--tangential boundary values where 
$\uX  \in D^- \longrightarrow \uXi \in \p D$ we find
$$
\lim_{\uX \rightarrow \uXi}\, g(\uX) =  0
$$
or
$$
-\,  \onehalf\,  h^n(\uXi) + \onehalf\,  \mH[h^n](\uXi) = 0
$$
again confirming that

$$
\mH[h^n](\uXi) = h^n(\uXi)
$$
}
\end{example}


\section{Quaternionic monogenicity}
\label{quaternionicmonogenicity}


A further refinement of hermitian monogenicity is obtained by taking the dimension to be a fourfold: $m=2n =4p$, reordering the variables as follows:
$$
(X_1,\ldots,X_{4p}) = (x_1,y_1,x_2,y_2,\ldots,x_{2p},y_{2p})
$$ 
and considering the hypercomplex structure $\mQ = (\mI_{4p},\mJ_{4p},\mK_{4p})$ on $\mR^{4p}$. This hypercomplex structure arises by introducing, next to the complex structure $\mI_{4p}$, a second one, $\mJ_{4p}$, given by
$$
\mJ_{4p} = \mbox{diag} \, \left ( \begin{array}{cccc} &  & 1 & \\ & & & -1 \\  -1 & & &  \\ & 1 & & \end{array} \right )
$$
Clearly $\mJ_{4p} \in \mbox{SO}(4p)$, with $\mJ_{4p}^2 = -E_{4p}$, and it anti--commutes with $\mI_{4p}$. 
A third $\mbox{SO}(4p)$--matrix 
$$
\mK_{4p} = \mI_{4p} \, \mJ_{4p} = - \mJ_{4p} \, \mI_{4p}
$$ 
then arises naturally,  for which $\mK_{4p}^2 = -E_{4p}$ and which anti--commutes with both $\mI_{4p}$ and $\mJ_{4p}$. Note that the representation of vectors is assumed to be by rows and the action of matrices on vectors thus is given by right multiplication, whence the above relation between the matrices $\mK$, $\mI$ and $\mJ$ in fact signifies that $\mK = \mJ \circ \mI$.\\[-2mm]

Next to the vector variable 
$$
\uX \ = \ \sum_{k=1}^n (x_k e_{2k-1} + y_k e_{2k}  + x_{k+1} e_{2k+1} + y_{k+1} e_{2k+2})
$$
we now introduce the rotated variables
\begin{eqnarray*}
\uX_{\mI} &=& \mI[\uX] = \sum_{k=1}^n (-y_k e_{2k-1} + x_k e_{2k}  - y_{k+1} e_{2k+1} + x_{k+1} e_{2k+2})\\
\uX_{\mJ} &=&  \mJ[\uX] = \sum_{k=1}^n (-x_{k+1} e_{2k-1} + y_{k+1} e_{2k} + x_k e_{2k+1} - y_k e_{2k+2} )\\
\uX_{\mK} &=& \mK[\uX] =  \sum_{k=1}^n (y_{k+1} e_{2k-1} + x_{k+1} e_{2k} - y_k e_{2k+1} - x_k e_{2k+2} )
\end{eqnarray*}
and we introduce the concept of quaternionic monogenicity by means of the following four operators: the Dirac operator
$$
\dirac \ = \ \sum_{k=1}^n (\partial_{x_k} e_{2k-1} + \partial_{y_k} e_{2k}  + \partial_{x_{k+1}} e_{2k+1} + \partial_{y_{k+1}} e_{2k+2})
$$
and the three additional rotated Dirac operators 
\begin{eqnarray*}
\partial_{\mI} &=& \mI_{4p}[\dirac] \ = \  \sum_{k=1}^n (-\partial_{y_k} e_{2k-1} + \partial_{x_k} e_{2k}  - \partial_{y_{k+1}} e_{2k+1} + \partial_{x_{k+1}} e_{2k+2})\\
\partial_{\mJ} &=& \mJ_{4p}[\dirac] \ = \  \sum_{k=1}^n (-\partial_{x_{k+1}} e_{2k-1} + \partial_{y_{k+1}} e_{2k} + \partial_{x_k} e_{2k+1} - \partial_{y_k} e_{2k+2} )\\
\partial_{\mK} &=& \mK_{4p}[\dirac] \ = \  \sum_{k=1}^n (\partial_{y_{k+1}} e_{2k-1} + \partial_{x_{k+1}} e_{2k} - \partial_{y_k} e_{2k+1} - \partial_{x_k} e_{2k+2} )
\end{eqnarray*}

\begin{definition}
\label{defqmonogenic}
A differentiable function $F: \mR^{4p} \longrightarrow \mS$ is called quaternionic monogenic in an open region $\Omega$ of $\mR^{4p}$, if and only if in that region $F$ is a solution of the system
$$
\{ \dirac F = 0, \dirac_\mI F = 0, \dirac_\mJ F = 0, \dirac_\mK F = 0 \}
$$
\end{definition}

Here too an alternative characterization is possible through complexification. In the actual dimension the hermitian vector variables read
\begin{eqnarray*}
\uz &=& - \frac{1}{2} ({\bf 1} - i \, \mI_{4p} ) [\uX]  \ = \ \sum_{j=1}^{p} ( z_{2j-1} \gf_{2j-1}  + z_{2j}  \gf_{2j} ) \\
\uzd &=& \phantom{-} \frac{1}{2} ({\bf 1} + i \, \mI_{4p} ) [\uX] \ = \  \sum_{j=1}^{p}\,  ( z_{2j-1}^c \gfd_{2j-1} + z_{2j}^c \gfd_{2j} )
\end{eqnarray*}
and their images under the action of $\mJ_{4p}$ turn out to be
\begin{eqnarray*}
\uzJ &= & \mJ_{4p}[\uz] \ = \  - \frac{1}{2} (\mJ_{4p} -  i \, \mK_{4p} ) [\uX] \ =\ \sum_{j=1}^p ( z_{2j} \gfd_{2j-1} -  z_{2j-1} \gfd_{2j} ) \\
\uzJd & = &  \mJ_{4p}[\uzd] \ = \  \phantom{-} \frac{1}{2} (\mJ_{4p} + i \, \mK_{4p} ) [\uX] \ = \  \sum_{j=1}^p ( z_{2j}^c \gf_{2j-1} -  z_{2j-1}^c \gf_{2j} )
\end{eqnarray*}
The corresponding quaternionic Dirac operators are
\begin{eqnarray*}
\upz & = & \phantom{-} \frac{1}{4} ( {\bf 1} + i \, \mI_{4p})[\dirac] \ = \ \sum_{j=1}^p ( \p_{z_{2j-1}} \gfd_{2j-1} + \p_{z_{2j}} \gfd_{2j} ), \\
\upzd & = & - \frac{1}{4} ( {\bf 1} - i \, \mI_{4p})[\dirac] \ = \ \sum_{j=1}^p ( \p_{z_{2j-1}^c} \gf_{2j-1} + \p_{z_{2j}^c} \gf_{2j} )\\
\upzJ & = & \phantom{-} \frac{1}{4} (\mJ_{4p} + i \, \mK_{4p} ) [\dirac] \ = \ \mJ_{4p} [\upz] \ = \ \sum_{j=1}^p ( \p_{z_{2j}} \gf_{2j-1} - \p_{z_{2j-1}} \gf_{2j} ), \\
\upzJd & = & - \frac{1}{4} (\mJ_{4p} - i \, \mK_{4p} ) [\dirac] \ = \ \mJ_{4p} [\upzd] \ = \ \sum_{j=1}^p ( \p_{z_{2j}^c} \gfd_{2j-1} - \p_{z_{2j-1}^c} \gfd_{2j} )
\end{eqnarray*}
For a function $F$ on $\mR^{4p} \cong \mC^{2n}$ the quaternionic monogenic system of Definition \ref{defqmonogenic} then is easily seen to be equivalent to the system
$$
\{ \upz F = 0, \upzd F = 0, \upzJ F = 0,  \upzJd F = 0 \}
$$
which can be shown to be invariant under the action of the symplectic group  Sp$(p)$. The basics of the quaternionic monogenic function theory were developed in \cite{paper1,paper2,paper3}. For group theoretical aspects we refer to \cite{vorig,mmas}.\\[-2mm]

In the real approach to quaternionic monogenicity we have the fundamental solutions
$$
E(\uX) = \frac{1}{a_{4p}}\, \frac{\overline{\uX}}{|\uX|^{4p}}, \quad E_\mI(\uX) = \frac{1}{a_{4p}}\, \frac{\overline{\uX}_\mI}{|\uX_\mI|^{4p}}, \quad E_\mJ(\uX) = \frac{1}{a_{4p}}\, \frac{\overline{\uX}_\mJ}{|\uX_\mJ|^{4p}}, \quad  E_\mK(\uX) = \frac{1}{a_{4p}}\, \frac{\overline{\uX}_\mK}{|\uX_\mK|^{4p}}
$$
for the operators $\dirac$, $\dirac_{\mI}$, $\dirac_{\mJ}$ and $\dirac_{\mK}$ respectively, where now $a_{4p}$ denotes the area of the unit sphere $S^{4p-1}$ in $\mR^{4p}$. By similar projections/decompositions as above they give rise to their quaternionic counterparts, explicitly given by:
\begin{eqnarray*}
E(\uz) &=&  -  E(\uX) + i\, E_\mI (\uX) \ = \  \frac{2}{a_{4p}}\, \frac{\uz}{\left| \uz \right|^{4p}}\\
E^\dagger(\uz) &=&  \phantom{- } E(\uX) + i\, E_\mI(\uX)  \ = \ \frac{2}{a_{4p}}\, \frac{\uz^\dagger}{\left| \uz \right|^{4p}}\\
E^J(\uz) &=& - E_\mJ(\uX) + i\, E_\mK(\uX) \ = \ \frac{2}{a_{4p}}\, \frac{\uz^J}{\left| \uz \right|^{4p}}\\
{E^\dagger}^J(\uz) &=& \phantom{-}  E_\mJ(\uX) + i\, E_\mK(\uX) \ = \ \frac{2}{a_{4p}}\, \frac{{\uz^\dagger}^J}{\left| \uz \right|^{4p}}
\end{eqnarray*}
However, as could be expected, the latter are not fundamental solutions for the quaternionic Dirac operators, whence, again, a circulant matrix approach has to be followed.  Recall that, in distributional sense, 
\begin{eqnarray*}
\upz E(\uz) &=& \frac{1}{2p} \beta\, \delta(\uz) + \frac{2}{a_{4p}} \beta\, \mbox{Fp} \frac{1}{|\uz|^{4p}} - \frac{2}{a_{4p}} (2p)\, \mbox{Fp} \frac{\uzd \uz}{|\uz|^{4p+2}}\\
\upzd E^\dagger(\uz) &=& \frac{1}{2p} (2p-\beta)\, \delta(\uz) + \frac{2}{a_{4p}} (2p-\beta)\, \mbox{Fp} \frac{1}{|\uz|^{4p}} - \frac{2}{a_{4p}} (2p)\, \mbox{Fp} \frac{\uz \uzd}{|\uz|^{4p+2}}
\end{eqnarray*}
whence
$$
\upz\, E(\uz) + \upzd\, E^\dagger(\uz) = \delta(\uz) 
$$
In addition we now have, again in distributional sense, that
\begin{eqnarray*}
\upzJ E^J(\uz) &=& \frac{1}{2p} (2p-\beta)\, \delta(\uz) + \frac{2}{a_{4p}} (2p-\beta)\, \mbox{Fp} \frac{1}{|\uz|^{4p}} - \frac{2}{a_{4p}} (2p)\, \mbox{Fp} \frac{\uzJd \uzJ}{|\uz|^{4p+2}}\\
\upzJd {E^\dagger}^J(\uz) &=& \frac{1}{2p} \beta\, \delta(\uz) + \frac{2}{a_{4p}} \beta\, \mbox{Fp} \frac{1}{|\uz|^{4p}} - \frac{2}{a_{4p}} (2p)\, \mbox{Fp} \frac{\uzJ \uzJd}{|\uz|^{4p+2}}
\end{eqnarray*}
whence
$$
\upzJ\, E^J(\uz) + \upzJd\, E^{\dagger J}(\uz) = \delta(\uz)
$$

The matrix  operator $\boldsymbol{ \mcD}$ given by
$$
\boldsymbol{\mathcal{D}} \ = \ \left( \begin{array}{llll}
\upz & \upzd & \upzJ & \upzJd \\[2mm] 
\upzJd & \upz & \upzd & \upzJ \\[2mm] 
\upzJ & \upzJd & \upz & \upzd \\[2mm]
\upzd & \upzJ & \upzJd & \upz
\end{array} \right)
$$
factorizes the matrix Laplace operator in the sense that 
$2 \boldsymbol{\mathcal{D}}  \boldsymbol{\mathcal{D}} ^\dagger = {\boldsymbol \Delta_{4p}}$. Introducing the matrices
$$
{\bf{E}}(\uz) \ = \  \left( \begin{array}{llll} E & E^\dagger & E^J & {E^\dagger}^J\\
{E^\dagger}^J & E & E^\dagger & E^J \\
E^J & {E^\dagger}^J & E & E^\dagger \\
E^\dagger & E^J & {E^\dagger}^J & E
\end{array} \right) \qquad \mbox{and} \qquad 
{\bd}(\uz) \ = \ \left( \begin{array}{llll} \delta(\uz) &0 & 0 & 0\\
0 & \delta(\uz) & 0 & 0\\
0 & 0 & \delta(\uz) & 0\\
0 & 0 & 0 & \delta(\uz) \\
\end{array} \right)
$$
it is easily obtained that
$$
\boldsymbol{\mathcal{D}}  \boldsymbol{E}^{\mbox{{\small T}}}(\uz) = 2\boldsymbol{\delta} (\uz) 
$$
whence a matrix fundamental solution has been found for the matrix Dirac operator $\boldsymbol{\mcD}$. Notice that in the hermitian case transposing the matrix ${\bf{E}}$ was not necessary, since a circulant $2 \times 2 $ matrix always is symmetric. A similar (yet slightly different) strategy was developed in \cite{ab}. \\[-2mm]

However, another approach is possible as well, since the actions $\upz E^\dagger(\uz)$, $\upzd E(\uz)$, $\upzJ {E^\dagger}^J(\uz)$ and $\upzJd E^J(\uz)$ all equal zero, implying that we can also consider
\begin{equation}
\boldsymbol{\mathcal{D}_{(\uz, \uzJ)}} \ = \ \left( \begin{array}{llll}
\upz & \upzd & 0 & 0\\[2mm] 
\upzd & \upz & 0 & 0 \\[2mm] 
0 & 0 & \upzJ & \upzJd \\[2mm]
0 & 0 & \upzJd & \upzJ
\end{array} \right) \qquad \mbox{and} \qquad 
\bf{E} \ = \  \left( \begin{array}{llll} E & E^\dagger & 0 & 0\\
E^\dagger & E & 0 & 0 \\
0 & 0 & E^J & {E^\dagger}^J \\
0 & 0 & {E^\dagger}^J & E^J
\end{array} \right)
\label{2ndattemptmat} 
\end{equation}
for which it holds that $4  \boldsymbol{\mathcal{D}_{(\uz, \uzJ)}} \boldsymbol{\mathcal{D}_{(\uz, \uzJ)}}^\dagger = {\boldsymbol \Delta}$ and 
\begin{equation}
\boldsymbol{\mathcal{D}_{(\uz, \uzJ)}} \boldsymbol{E}(\uz) = \boldsymbol{\delta} (\uz)
\label{2ndattemptsol}
\end{equation}

In the next section we will find out which of both possibilities is best suited for establishing a Cauchy Integral Formula in the quaternionic Clifford setting.


\section{The Cauchy Integral Formula in the quaternionic case}
\label{cauchyquaternionic}


In order to make a deliberate choice between both approaches sketched in the preceding section, we will first have a look at the underlying group symmetry for the quaternionic monogenic function theory. To that end let us again consider functions taking values in complex spinor space, which now is given by
$$
\mS = \mC_{4p} \, I \cong \mC_{2p} \, I
$$
and which is realized by means of the primitive idempotent $I = I_1 \ldots I_{2p}$, with $I_j = \gf_j \gf_j^{\dagger}$, $j=1,\ldots,2p$. We already know that, as a U$(n)$--module, complex spinor space  decomposes into homogeneous parts as
$$
\mS = \bigoplus_{r=0}^{2p} \mS^r = \bigoplus_{r=0}^{2p} (\mC\Lambda_{2p}^\dagger)^{(r)} I
$$
An important observation is that the spaces $\mS^r$, which are invariant and irreducible U$(n)$ modules, are reducible under the action of the fundamental symmetry group Sp$(p)$.  \\[-2mm]

It still holds that a spinor valued function $g$ is quaternionic monogenic if and only if all its homogeneous parts $g^{r}$ are. However for a fixed component $g^{r}$ quaternionic monogenicity is not equivalent to monogenicity. Yet we have the following result, see \cite{paper2}.

\begin{lemma}
For a function $g^r$ defined on (a domain in) $\mR^{4p} \cong \mC^{2n}$ and taking values in $\mS^r$, $r \in \{0,1,\ldots,2p\}$, it holds that $g^r$ is quaternionic monogenic if and only if it is simultaneously $\dirac$ and $\dirac_{\mJ}$ monogenic.
\end{lemma}

This result shows that, in view of establishing a Cauchy Integral Formula,  the second attempt (\ref{2ndattemptmat})--(\ref{2ndattemptsol}) is the right one to pursue, since the structure of the involved matrices reflects the importance of $\dirac$ and $\dirac_{\mJ}$ monogenicity in this setting. We thus introduce the concept of matricial quaternionic monogenicity as follows.
\begin{definition}
A block diagonal matrix function
$$
\bG \ = \  \left( \begin{array}{cccc} g_1 & g_2 & 0 & 0 \\ g_2 & g_1 & 0 & 0 \\ 0 & 0 & g_3 & g_4 \\ 0 & 0 & g_4 & g_3  \end{array} \right) 
$$
with continuously differentiable entries $g_1,g_2,g_3,g_4$ in (a domain in) $\mR^{4p} \cong \mC^{2p}$ taking values in (subspaces of) $\mC_{4p}$  is called quaternionic monogenic if and only if it satisfies the system 
$$
\boldsymbol{\mathcal{D}_{(\uz, \uzJ)}}  \bG = \bO
$$ 
\end{definition}

Note that in the case of a diagonal matrix function $\bG_0$ with $g_1=g=g_3$ and $g_2=0=g_4$, the quaternionic monogenicity of $\bG_0$ coincides with the quaternionic monogenicity of $g$, which is not the case in general.\\[-2mm]

Next we define, in a similar way, the matrix differential form
$$
\bdS_{(\uz,\uzJ)} \ = \ \left( \begin{array}{llll}  d\sigma_{\uz} &  d\sigma_{\uz^{\dagger}} & 0 & 0 \\ d\sigma_{\uz^{\dagger}}  & d\sigma_{\uz} & 0 & 0 \\
0 & 0 & d\sigma_{\uz^J} & d\sigma_{{\uz^\dagger}^J}\\
0 & 0 & d\sigma_{{\uz^\dagger}^J} & d\sigma_{\uz^J}  \\
\end{array} \right)
$$
Here, as above, we have introduced
\begin{eqnarray*}
d\sigma_{\uz^{\phantom{\dagger}}}  &=&  \phantom{-} \sum_{j=1}^{p} \left ( \gf_{2j-1}^\dagger \, \widehat{dz_{2j-1}} + \gf_{2j}^\dagger \, \widehat{dz_{2j}} \right ) \\
d\sigma_{\uz^{\dagger}}  &=&  - \sum_{j=1}^{p} \left ( \gf_{2j-1} \,  \widehat{dz_{2j-1}^c}  + \gf_{2j} \, \widehat{dz_{2j}^c} \right )
\end{eqnarray*}
where the notations $\widehat{dz_k}$ and $\widehat{dz_k^c}$ keep their original definition, see (\ref{hatdz})--(\ref{hatdzdagger}), whence
\begin{eqnarray*}
d\sigma_{\uz^{\phantom{\dagger}}}  & = & (-i)^{2p} 2^{2p-1} \pi^-[d\sigma_{\uX}] \ = \  - \frac{1}{2}  (-i)^{2p} 2^{2p-1}  \left ( d\sigma_{\uX} - i \, d\sigma_{\uX_{\mI}} \right ) \\
d\sigma_{\uz^{\dagger}}  & = & (-i)^{2p} 2^{p-1} \pi^+[d\sigma_{\uX}] \ = \   \phantom{-} \frac{1}{2}  (-i)^{2p} 2^{2p-1}  \left ( d\sigma_{\uX} + i \, d\sigma_{\uX_{\mI}} \right ) 
\end{eqnarray*}
and in a similar way we have defined
\begin{eqnarray*}
d\sigma_{{\uzJ}^{\phantom{\dagger}}}  &=&  \phantom{-} \sum_{j=1}^{p} \left ( \gf_{2j-1} \, \widehat{dz_{2j}} - \gf_{2j} \, \widehat{dz_{2j-1}} \right ) \\
d\sigma_{\uzJd}  &=&  - \sum_{j=1}^{p} \left ( \gf_{2j-1}^\dagger \,  \widehat{dz_{2j}^c}  - \gf_{2j}^\dagger \, \widehat{dz_{2j-1}^c} \right )
\end{eqnarray*}
or, expressed in the original real variables
\begin{eqnarray*}
d\sigma_{\uz^{\phantom{\dagger}}}  & = & J[d\sigma_{\uz}] \ = \  - \frac{1}{2}  (-i)^{2p} 2^{2p-1}  \left ( d\sigma_{\uX_{\mJ}} - i \, d\sigma_{\uX_{\mK}} \right ) \\
d\sigma_{\uz^{\dagger}}  & = &J[d\sigma_{\uzd}] \ = \   \phantom{-} \frac{1}{2}  (-i)^{2p} 2^{2p-1}  \left ( d\sigma_{\uX_{\mJ}} + i \, d\sigma_{\uX_{\mK}} \right ) 
\end{eqnarray*}

The Cauchy Integral Formula for matrix quaternionic monogenicity readily follows; its proof is similar to the one in the hermitian monogenic setting.
\begin{proposition}
Let the block diagonal matrix function $\bG$ be quaternionic monogenic in an open neighbourhood of $\overline{D}$, $D \subset \mR^{4p}$ being a bounded domain in $\mR^{4p}$ with smooth boundary $\partial D$. Then it holds that
$$
 \bG(\uX) = \frac{1}{(-4)^{p}} \int_{\partial D} \boldsymbol{E} (\uv - \uz)\, \bdS_{(\uv,\uvJ)} \, \bG(\uY) \quad ,  \quad \uX \in \overset{\circ}{D} 
$$
where $\uv$ is the hermitian vector variable corresponding to $\uY \in \partial D$, and $\uz$ is the one corresponding to $\uX \in \overset{\circ}{D}$. 
\end{proposition}

Again note that the multiplicative constant appearing at the right hand side originates from the re-ordering of $4p$ real variables into $2p$ complex planes. \\[-2mm]

Taking for $\bG$ the diagonal matrix function $\bG_0$, the above result reduces to a genuine Cauchy Integral Formula for the quaternionic monogenic function $g$.
\begin{corollary}
\label{corCIFquat}
Let the function $g$ be quaternionic monogenic in an open neighbourhood of $\overline{D}$, $D \subset \mR^{4p}$ being a bounded domain in $\mR^{4p}$ with smooth boundary $\partial D$. Then one has the two reproducing formul{\ae}
\begin{eqnarray}
\label{CIFquat1} g(\uX) &=& \frac{1}{(-4)^p} \int_{\partial D}  \left [ E(\uv-\uz)\,  d\sigma_{\uv} + E^\dagger(\uv-\uz)\,  d\sigma_{\uv^{\dagger}} \right ] \, g(\uY) \quad ,  \quad \uX \in \overset{\circ}{D}  \\
\label{CIFquat2} g(\uX) &=& \frac{1}{(-4)^p} \int_{\partial D}  \left [ E^J(\uv-\uz)\,  d\sigma_{\uvJ} + {E^\dagger}^J(\uv-\uz)\,  d\sigma_{\uvJd} \right ] \, g(\uY)  \quad ,  \quad \uX \in \overset{\circ}{D} 
\end{eqnarray}
and the two additional integral identities
\begin{eqnarray}
\label{IDquat1} \int_{\partial D}  \left [ E(\uv-\uz)\,  d\sigma_{\uv^\dagger} + E^\dagger(\uv-\uz)\,  d\sigma_{\uv} \right ] \, g(\uY) & = & 0 \quad ,  \quad \uX \in \overset{\circ}{D}\\
\label{IDquat2} \int_{\partial D}  \left [ E^J(\uv-\uz)\,  d\sigma_{\uvJd} + {E^\dagger}^J(\uv-\uz)\,  d\sigma_{\uvJ} \right ] \, g(\uY) & = & 0 \quad ,  \quad \uX \in \overset{\circ}{D}
\end{eqnarray}
\end{corollary} 

An alternative proof of Corollary \ref{corCIFquat} is obtained, as has been done explicitly in the hermitian case, by splitting a spinor valued function in its homogeneous components, writing down, for each component, the Cauchy Integral Formul{\ae} for $\partial$ and $\partial_{\mJ}$ monogenicity, while invoking the structural decompositions for all building blocks involved and the subsequent splitting of the values, and finally adding up all intermediate results. In this way the reproducing formul{\ae} (\ref{CIFquat1}) and (\ref{CIFquat2}) are recovered, while the integral identities (\ref{IDquat1}) an (\ref{IDquat2}) are replaced by the stronger results:
\begin{align*}
\int_{\partial D}   E(\uv-\uz) d\sigma_{\uv^\dagger} \, g(\uY)  \ \ = \ \ &0 \ \ = \ \ E^\dagger(\uv-\uz) d\sigma_{\uv}  \, g(\uY)  \quad ,  \quad \uX \in \overset{\circ}{D}\\
\int_{\partial D}   E^J(\uv-\uz) d\sigma_{\uvJd}  \, g(\uY)  \ \ = \ \   &0 \ \ = \ \  {E^\dagger}^J(\uv-\uz) d\sigma_{\uvJ} \, g(\uY)  \quad ,  \quad \uX \in \overset{\circ}{D}
\end{align*}\hfill\\

Now recall that interesting results were obtained in the hermitian framework by restricting the values of the considered functions to the different homogenous parts of spinor space, which are suggested by the U$(n)$ symmetry. For quaternionic monogenics the underlying Sp$(p)$ invariance has not yet been fully exploited, since the homogeneous parts of spinor space are reducible under Sp$(p)$ and split further into so-called symplectic cells, see e.g.\ \cite{paper1,paper3}. This splitting is caused by the action of the multiplication operators
$$
P=\gf_2 \gf_1 + \gf_4 \gf_3 + \ldots + \gf_{2p} \gf_{2p-1}, \qquad Q = \gfd_1 \gfd_2 + \gfd_3 \gfd_4 + \ldots + \gfd_{2p-1} \gfd_{2p}
$$
for which we define, for $r=0,\ldots,p$, the kernel spaces
$$
\mS_r^r = \mbox{Ker} \, P |_{\mS^r}, \qquad \mS_r^{2p-r} = \mbox{Ker} \, Q |_{\mS^{2p-r}}
$$
and for $k=0,\ldots,p-r$, the subspaces obtained by iterative action of $Q$ on the kernel of $P$ and vice versa:
$$
\mS_r^{r+2k} = Q^k \, \mS^r_r, \qquad \mS_r^{2p-r-2k} = P^k \, \mS^{2p-r}_r
$$
It was shown in \cite{eel2, paper1} that, for all $r=0,\ldots,p$,
$$
\mS^r =\bigoplus_{j=0}^{\lfloor \frac{r}{2} \rfloor} \mS_{r-2j}^r, \qquad
\mS^{2p-r} =\bigoplus_{j=0}^{\lfloor \frac{r}{2} \rfloor} \mS_{r-2j}^{2p-r}
$$
and each of the symplectic cells $\mS_s^r$ in the above decompositions is an irreducible $\mbox{Sp}(p)$--representation. Whence we can now decompose a function $F : \mR^{4p} \longrightarrow \mS$ into components taking values in these symplectic cells:
$$
F = \sum_{r=0}^n \, F^r =  \sum_{r=0}^n \, \sum_s F^r_s, \qquad F^r_s : \mR^{4p} \longrightarrow \mS^r_s
$$
The quaternionic monogenicity of $F$ then is shown to be equivalent with the quaternionic monogenicity of each of its components $F^r_s$, entailing an even further refinement of the results obtained above.\\[-2mm]


\section{The $\gosp(4|2)$--monogenic framework}
\label{osp42}


In \cite{mmas} it was shown that, from a group theoretical point of view, the definition of quaternionic monogenicity is not the best possible one. For instance, spaces $\mcQ^{r,s}_{a,b}$ of quaternionic monogenic bi--homogeneous polynomials of bi--degree $(a,b)$ with values in the symplectic cell $\mS^r_s$, still remain reducible under the action of the symplectic group Sp$(p)$, an unfortunate situation. This has lead to the definition of so-called $\gosp(4|2)$--monogenicity in \cite{vorig,mmas}, where a function, apart from being quaternionic monogenic, is requested to be in the kernel of the above mentioned multiplication operator $P$:
$$
P = \gf_2 \gf_1 + \cdots + \gf_{2p} \gf_{2p-1}
$$ 
and in the kernel of the Euler like scalar differential operator
$$
\SE =  \sum_{k=1}^p  \, z_{2k-1} \, \p_{z_{2k}^c} - z_{2k} \, \p_{z_{2k-1}^c}
$$
which arises when computing the anti-commutators of all operators in the odd part of the involved Lie superalgebra $\gosp(4|2)$.

\begin{definition}
A function $f$ is $\gosp(4|2)$--monogenic in an open region $\Omega$ of $\mR^{4p}$ if in $\Omega$ it belongs to the kernel of the six operators: $\upz$, $\upzd$, $\upzJ$, $\upzJd$, $P$ and $\SE$.
\end{definition}

\noindent
As Ker$P = \bigoplus_{r=1}^p\, \mS_r^r$, an $\gosp(4|2)$--monogenic function $f$ ought to take its values in $\bigoplus_{r=1}^p\, \mS_r^r$, whence it can be decomposed as
$$
f = \sum_{r=1}^p\, f_r^r
$$
and such a function is quaternionic monogenic if and only if all the components $f_r^r$ are. As for each $r$ it trivially holds that $\upzd\, f_r^r = \upzJd\, f_r^r = 0$, the quaternionic monogenicity of $f$ is equivalent with the system \{$\upz\, f_r^r = 0, \upzJ\, f_r^r = 0 $\}.\\

Exploiting the results on quaternionic monogenic functions, we are able to establish a Cauchy Integral Formula for $\gosp(4|2)$--monogenic functions.

\begin{theorem}
\label{CIFosp}
Let the function $g$ be $\gosp(4|2)$--monogenic in an open neighbourhood $\Omega$ of $\overline{D}$, $D$ being a bounded domain $D \subset \mR^{4p}$ with smooth boundary $\p D$.  Then the following representation formulae hold:
\begin{align}
\label{osprep1} g(\uX) &=     \frac{1}{(-2i)^n}\,  \int_{\p D}\, \left(  E(\uv-\uz)\, d\sigma_{\uv} +    E^{\dagger}(\uv-\uz)\, d\sigma_{\uvd} \right)\, g(\uY) \quad , \quad  \uX \in \overset{\circ}{D}\\
\label{osprep2} g(\uX) &=     \frac{1}{(-2i)^n}\,  \int_{\p D}\, \left(  E^J(\uv-\uz)\, d\sigma_{\uvJ} +    E^{J \dagger}(\uv-\uz)\, d\sigma_{\uvJd} \right)\, g(\uY) \quad , \quad  \uX \in \overset{\circ}{D}
\end{align}
together with the integral identities
\begin{align}
\label{ospid1} 0 &=   \int_{\p D}\,   E^{\dagger}(\uv-\uz)\, d\sigma_{\uv} \, g(\uY) \quad , \quad  \uX \in \overset{\circ}{D}\\
\label{ospid2} 0 &=   \int_{\p D}\,   E^{J}(\uv-\uz)\, d\sigma_{\uvJd} \, g(\uY) \quad , \quad  \uX \in \overset{\circ}{D}\\
\label{ospid3} 0 &=   \int_{\p D}\,  \left( \SE\,  E(\uv-\uz)\, d\sigma_{\uv} +    \SE\, E^{\dagger}(\uv-\uz)\, d\sigma_{\uvd} \right)\, g(\uY) \quad , \quad  \uX \in \overset{\circ}{D}\\
\label{ospid4} 0 &=   \int_{\p D}\,  \left( \SE\,  E^J(\uv-\uz)\, d\sigma_{\uvJ} +    \SE\, E^{J \dagger}(\uv-\uz)\, d\sigma_{\uvJd} \right)\, g(\uY) \quad , \quad  \uX \in \overset{\circ}{D}
\end{align}
\end{theorem}

\pf
(i) The representation formulae (\ref{osprep1}) and (\ref{osprep2}) are due to the quaternionic monogenicity of the function $g$ in $\Omega$.\\
(ii) The identities (\ref{ospid1}) and (\ref{ospid2}) hold for the same reason, taking into account that, trivially,
$$
0 =   \int_{\p D}\,   E(\uv-\uz)\, d\sigma_{\uvd} \, g(\uY) \quad , \quad  \uX \in \overset{\circ}{D}
$$
and
$$
 0 =   \int_{\p D}\,   E^{J \dagger}(\uv-\uz)\, d\sigma_{\uvJ} \, g(\uY) \quad , \quad  \uX \in \overset{\circ}{D}
$$
(iii) The identities (\ref{ospid3}) and (\ref{ospid4}) are the result of the action of the operator $\SE$ on the representation formulae (\ref{osprep1}) and (\ref{osprep2}), taking into account that $\SE g = 0$.
\qed

\begin{remark}
{\rm
The kernel functions appearing in the identities (\ref{ospid3}) and (\ref{ospid4}) may be calculated explicitly. We obtain:
\begin{align*}
\SE\, E(\uv-\uz) &=   \frac{4p}{a_{4p}}\, \frac{\uv-\uz}{\rho^{4p+2}}\, \sum_{j=1}^p\, (z_{2j-1}v_{2j} - z_{2j}v_{2j-1}) \\
\SE\, E^{\dagger}(\uv-\uz) &=    \frac{4p}{a_{4p}}\, \frac{\uvd-\uzd}{\rho^{4p+2}}\, \sum_{j=1}^p\, (z_{2j-1}v_{2j} - z_{2j}v_{2j-1}) + \frac{2}{a_{4p}}\, \frac{1}{\rho^{4p}}\, \uzJ\\
\SE\, E^J(\uv-\uz) &=    \frac{4p}{a_{4p}}\, \frac{\uvJ-\uzJ}{\rho^{4p+2}}\, \sum_{j=1}^p\, (z_{2j-1}v_{2j} - z_{2j}v_{2j-1})  \\
\SE\, E^{J \dagger}(\uv-\uz) &=  \frac{4p}{a_{4p}}\, \frac{\uvJd-\uzJd}{\rho^{4p+2}}\, \sum_{j=1}^p\, (z_{2j-1}v_{2j} - z_{2j}v_{2j-1}) - \frac{2}{a_{4p}}\, \frac{1}{\rho^{4p}}\, \uz
\end{align*}
}
\end{remark}

Now we try to establish the concept of Cauchy transform in the context of $\gosp(4|2)$--monogenicity. As usual we start with a bounded domain $D$ in $\mR^{4p}$ with a smooth boundary $\p D$ and a continuous function $h$ on $\p D$, which, from the start, is assumed to take values in Ker$P$:
$$
h = \sum_{r=0}^p\, h_r^r \quad , \quad h_r^r : \p D \rightarrow \mS^r_r
$$
We define the following four functions in $\overset{\circ}{D}$:

\begin{align*}
g_1(\uX) &= \int_{\p D}\, \left( E(\uv-\uz)\, d \sigma_{\uv} + E^{\dagger}(\uv-\uz)\, d \sigma_{\uvd}   \right)\, h(\uY) \\
g_2(\uX) &= \int_{\p D}\,  E^{\dagger}(\uv-\uz)\, d \sigma_{\uv}\, h(\uY) \\
g_3(\uX) &= \int_{\p D}\, \left( E^J(\uv-\uz)\, d \sigma_{\uvJ} + E^{J \dagger}(\uv-\uz)\, d \sigma_{\uvJd}   \right)\, h(\uY) \\
g_4(\uX) &= \int_{\p D}\,  E^J(\uv-\uz)\, d \sigma_{\uvJd}\, h(\uY)
\end{align*}
Clearly $g_2$ and $g_4$ take values in $\bigoplus_{r=0}^p\, \mS^{r+2}_r$, whence we make the first assumption that the boundary function $h$ satisfies the conditions: $g_2 = 0$ and $g_4 = 0, \forall \uX \in D^+$. The functions $g_1$ and $g_3$ take values in Ker$P$. In order that $g_1$ and $g_3$ should belong to Ker $\mcE$ we make the second assumption that  the boundary function $h$ satisfies the conditions:
$$
\int_{\p D}\, \left( \mcE\, E(\uv-\uz)\, d \sigma_{\uv} + \mcE\, E^{\dagger}(\uv-\uz)\, d \sigma_{\uvd}   \right)\, h(\uY)  = 0 \quad , \quad  \uX \in D^+
$$
and
$$
\int_{\p D}\, \left( \mcE\, E^J(\uv-\uz)\, d \sigma_{\uvJ} + \mcE\, E^{J \dagger}(\uv-\uz)\, d \sigma_{\uvJd}   \right)\, h(\uY) = 0 \quad , \quad  \uX \in D^+
$$
Further it holds that, in a trivial way, $\upzd\, g_1 = \upzJ\, g_1 = 0$ and $\upzd\, g_3 = \upzJ\, g_3 = 0$. Moreover, writing the above expressions for $g_1, g_2=0, g_3, g_4=0$ in a matricial form and acting with the matricial Dirac operator $\boldsymbol{\mathcal{D}}$  leads to $\upz\, g_1 = 0$ and $\upzJd\, g_3 = 0$. So, under the first and second assumption, it holds  that:\\
(i) if it would be so that $g_1 = g_3$, then this function is $\gosp(4|2)$--monogenic  in $D^+$;\\
(ii) if $\upzJd\, g_1 = 0$ then $g_1$ is $\gosp(4|2)$--monogenic  in $D^+$;\\
(iii) if  $\upz\, g_3 = 0$, then $g_3$ is $\gosp(4|2)$--monogenic  in $D^+$.\\[2mm]
The conclusion is that the construction of an $\gosp(4|2)$--monogenic function in the interior of a bounded domain via the Cauchy transform of a continuous boundary function is not straightforward. But that it is indeed possible is illustrated by the following example.

\begin{example}
{\rm
Take a scalar polynomial $V(z_1, z_1^c, \ldots, z_n, z_n^c)$ in Ker\,$\mcE$; it is even possible to take a harmonic such polynomial, see \cite{paper3}. Choose a vector $f$  in $\mS_p^p$ and put
$H_p^p = V f$. Then $H_p^p$ is quaternionic monogenic (trivially), in Ker$P$ (trivially) and in Ker $\mcE$ due to the choice of $V$; in other words: $H_p^p$ is $\gosp(4|2)$--monogenic in $\mR^{4p}$. It follows that the identities (\ref{ospid3}) and (\ref{ospid4}) hold for $H_p^p$, which, putting $h_p^p = H_p^p|_{\p D}$ read:
\begin{align*}
0 &=   \int_{\p D}\,  \left( \SE\,  E(\uv-\uz)\, d\sigma_{\uv} +    \SE\, E^{\dagger}(\uv-\uz)\, d\sigma_{\uvd} \right)\, h_p^p(\uY)\\
0 &=   \int_{\p D}\,  \left( \SE\,  E^J(\uv-\uz)\, d\sigma_{\uvJ} +    \SE\, E^{J \dagger}(\uv-\uz)\, d\sigma_{\uvJd} \right)\, h_p^p(\uY)
\end{align*}
showing that $h_p^p$ satisfies the second assumption.\\
Now construct the above mentioned four functions $g_1, g_2, g_3$ and $g_4$ with $h_p^p$ as the boundary function. In a trivial way it holds that $g_2 = g_4 = 0$, which means that also the first assumption is satisfied for the boundary function $h_p^p$.\\
From the general theory we know that 
\begin{align*}
\upz\, g_1 = 0 \qquad \upzd\, g_1 &= 0 \qquad \upzJ\, g_1 = 0\\
\upzd\, g_3 &= 0 \qquad \upzJ\, g_3 = 0 \qquad \upzJd\, g_3 = 0
\end{align*}
But for this boundary function $h_p^p$ it holds trivially that $\upzJd\, g_1 = 0$ and $\upz\, g_3 = 0$, whence 
$$
g_1(\uX) = \int_{\p D}\, \left( E(\uv-\uz)\, d \sigma_{\uv} + E^{\dagger}(\uv-\uz)\, d \sigma_{\uvd}   \right)\, V(\uY) f
$$
 and 
$$
g_3(\uX) = \int_{\p D}\, \left( E^J(\uv-\uz)\, d \sigma_{\uvJ} + E^{J \dagger}(\uv-\uz)\, d \sigma_{\uvJd}   \right)\, V(\uY) f
$$
 both are $\gosp(4|2)$--monogenic in $D^+$.
}
\end{example}


\section*{Acknowledgement}


R.\ L\'{a}vi\v{c}ka and V.\ Sou\v{c}ek gratefully acknowledge support by the Czech Grant Agency through grant GA CR 17-01171S. \\
This paper was partly written during a scientific stay of R.\ L\'{a}vi\v{c}ka at the former Department of Mathematical Analysis of Ghent University. R.\ L\'{a}vi\v{c}ka expresses his gratitude for the generous support and great hospitality during his stay.


\end{document}